\documentclass[11pt,a4paper]{article}

\usepackage[centertags,leqno]{amsmath}
\usepackage{amsthm}
\usepackage{amssymb,amsfonts,stmaryrd}
\usepackage{graphicx}
\usepackage{hyperref}

\allowdisplaybreaks

\numberwithin{equation}{section}

\newtheorem{thm}{Theorem}[section]
\newtheorem{prop}[thm]{Proposition}

\newtheorem{lemma}[thm]{Lemma}
\newtheorem{cor}[thm]{Corollary}
\newtheorem{Def}[thm]{Definition}
\theoremstyle{definition} 
\newtheorem{rem}[thm]{Remark}

\newenvironment{pf}{\textit{Proof.}}{\hfill$\boxbox$\\}
\newenvironment{pf-lb}{\textit{Proof.}}{\hfill$\boxbox$}

\newcommand{\R}{\mathbb{R}}

\newcommand{\Z}{\mathbb{Z}}
\newcommand{\C}{\mathbb{C}}
\renewcommand{\H}{\mathbb{H}}

\newcommand{\E}{\mathcal{E}}

\newcommand{\Isom}{\operatorname{Isom^+}}

\renewcommand{\sl}{\mathfrak{sl}}

\newcommand{\SL}{\operatorname{SL}}

\newcommand{\SU}{\operatorname{SU}}
\newcommand{\PSL}{\operatorname{PSL}}

\newcommand{\SO}{\operatorname{SO}}
\newcommand{\Ad}{\operatorname{Ad}}

\newcommand{\hol}{\operatorname{hol}}
\newcommand{\id}{\operatorname{id}}

\newcommand{\tr}{\operatorname{tr}}
\newcommand{\Tr}{\operatorname{Tr}}

\renewcommand{\Im}{\operatorname{Im}}
\renewcommand{\Re}{\operatorname{Re}}

\newcommand{\Defo}{\operatorname{Def}}

\begin{document}
\title{Complex twist flows on surface group representations and the local shape of the deformation space of hyperbolic cone-$3$-manifolds}

\author{Gr\'egoire Montcouquiol \footnote{Univ Paris-Sud, Laboratoire de Math\'ematiques, UMR8628, Orsay, F-91405; CNRS, Orsay, F 91405. Partially supported by the ANR programs GeomEinstein 06-BLAN-0154 and RepSurf 06-BLAN-0311} \and Hartmut Wei{\ss} \footnote{LMU M\"unchen, Mathematisches Institut, Theresienstr. 39, 80333 M\"unchen. Partially supported by the DFG Schwerpunktprogramm 1154 ``Globale Differentialgeometrie''}}

\maketitle

\begin{abstract}
In the former articles \cite{Montcouq2,We2}, it was independently proven by the authors that the space of hyperbolic cone-$3$-manifolds with cone angles less than $2\pi$ and fixed singular locus is locally parametrized by the cone angles. In this sequel, we investigate the local shape of the deformation space when the singular locus is no longer fixed, i.e.~when the singular vertices can be split. We show that the different possible splittings correspond to specific pair-of-pants decompositions of the smooth parts of the links of the singular vertices, and that under suitable assumptions the corresponding subspace of deformations is parametrized by the cone angles of the original edges and the lengths of the new ones.

\medskip
\noindent \emph{AMS classification:} 57M50, 58D27, 53C35
\end{abstract}


\section{Introduction}

It is well-known, since the fundamental works of Weil and Mostow \cite{Mostow,Weil}, that closed hyperbolic $3$-manifolds are rigid; the non-compact, complete case is also well understood. However for incomplete hyperbolic metrics, the situation is more complicated, and in this setting, it is natural to look at the metric completion of the manifold. For simplicity, we will consider the case of a hyperbolic metric $g$ on a $3$-manifold $M$ which is the interior of a compact manifold with boundary $\bar M$. Then there are, broadly speaking, two distinct types of situations depending on whether or not $g$ extends to a metric on the boundary. In the latter, the metric degenerates on $\partial \bar M$: this can happen in any number of ways, but cone-manifolds provide arguably the simplest and most interesting class of examples. In this setting, the boundary $\partial \bar M$ collapses to a (possibly disconnected) geodesic graph, called the \emph{singular locus} $\Sigma$, along which the metric exhibits a simple, ``cone-like'' singularity. The completion of $M$ is then a length space $X$, without boundary, such that $X\setminus \Sigma = (M,g)$; for more details and precise definitions, see \cite{BLP,CHK,Maz-GM,ThurstonShapes}. Cone-manifolds are particularly interesting to study, since they arise naturally in many different contexts:  they first appeared as deformations of complete, cusped hyperbolic $3$-manifolds \cite{ThurstonGeom}; they are the natural models for orbifolds, and the theory of their deformations plays a prominent role in the proof of the orbifold version of the geometrization theorem \cite{BLP,CHK}; as doubles of polyhedra, they form a natural framework for the resolution of Stoker's problem \cite{Maz-GM,Stoker}; finally, let us mention their use as models for space-times with massive point-like particles \cite{Kra-Sch}.

\medskip

Throughout this article, $X$ will denote a closed, orientable hyperbolic cone-$3$-manifold, and $\Sigma$ its singular locus. The homeomorphism type of the pair $(X,\Sigma)$ is called the {\em topological type} of the cone-manifold. The {\em smooth} or {\em regular part} of $X$ is the (incomplete) hyperbolic manifold $M=X \setminus \Sigma$. To each edge $e_i$ of the singular locus, one can assign a quantity called its \emph{cone angle}: it is defined as the positive real number $\alpha_i$ such that in cylindrical coordinates, the metric near any point of the interior of $e_i$ is expressed as 
$$g=dr^2 + \sinh^2(r)\,d\theta^2 + \cosh^2(r)\,dz^2, \quad r\in (0,\epsilon),\, \theta \in \R/\alpha_i \Z, \, z \in (-\epsilon,\epsilon).$$
If all cone angles are less than $2\pi$, then $X$ is a metric space with curvature bounded below by $-1$ in the triangle comparison sense. Furthermore, if all cone angles are less than or equal to $\pi$, then the vertices of $\Sigma$ are at most trivalent, which simplifies significantly the deformation theory. In general, the valence of the singular vertices can be arbitrarily high.

We are interested in the space of hyperbolic cone-manifold structures in the case where the cone angles are less than $2\pi$, i.e.~contained in the interval $(0,2\pi)$. More precisely, let $C_{-1}(X,\Sigma)$ denote the space of hyperbolic cone-manifold structures on $X$ of fixed topological type $(X,\Sigma)$; it is topologized as a subspace of the deformation space $\Defo(M)$ of incomplete hyperbolic structures on $M$. Let $N$ be the number of edges in $\Sigma$. In their seminal article \cite{HoKe}, C.D.~Hodgson and S.P.~Kerckhoff showed that if the singular locus is a link, i.e.~it does {\em not} contain vertices, then the map
$$
\alpha=(\alpha_1, \ldots, \alpha_N):C_{-1}(X,\Sigma) \rightarrow (0,2\pi)^N
$$
sending a hyperbolic cone-manifold structure to the vector of its cone angles is a local homeomorphism at the given structure. Hodgson and Kerckhoff proved more generally that $\Defo(M)$ is locally a smooth manifold at the given structure and that $\dim_\R \Defo(M)=2N$; cone-manifold structures are then identified to lie on a half-dimensional submanifold. In \cite{We1}, the second author showed that the same is true when vertices are allowed, if the cone angles are less than or equal to $\pi$. The general case has been open until recently, when the following has been proven independently in \cite{Montcouq2} and in \cite{We2}:

\begin{thm}\label{thm:former}
Let $X$ be a hyperbolic cone-$3$-manifold with cone angles less than $2\pi$. Then the map 
$$
\alpha=(\alpha_1, \ldots, \alpha_N): C_{-1}(X,\Sigma) \rightarrow (0,2\pi)^N
$$ 
is a local homeomorphism at the given structure.
\end{thm}

A similar result applies in the spherical and the Euclidean cases, at least when the cone angles are less than or equal to $\pi$, cf.~\cite{PoWe,We1}. However, if the cone angle bound is $2\pi$, then there exist counter-examples in the spherical case \cite{Porti,SchlenkerPoly}, whereas a weaker infinitesimal rigidity result still holds in the Euclidean case \cite{Maz-GM}.

The techniques used to show Theorem \ref{thm:former} in \cite{Montcouq2} and \cite{We2} are somewhat different: in the former, the starting point is an infinitesimal rigidity theorem for Einstein deformations obtained by the first author in a joint work with R.~Mazzeo, cf.~\cite{Maz-GM}, while the latter proceeds along the lines of \cite{HoKe,We1}. There the main technical ingredient is a vanishing theorem for $L^2$-cohomology with values in a certain flat bundle, see section \ref{sec:L2}. Common to both approaches is the use of the variety of representations $\rho: \pi_1 M \rightarrow \SL_2(\C)$ to pass from an infinitesimal rigidity statement to a local deformation theorem. Indeed, we recall that the map that sends a hyperbolic structure to its holonomy representation $\hol: \pi_1 M \rightarrow \SL_2(\C)$ induces a local homeomorphism 
\begin{equation}\label{eq:hol}
[\hol] : \Defo(M) \to X(\pi_1M,\SL_2(\C))
\end{equation}
where the right-hand side is the space of representations $\rho: \pi_1M \rightarrow \SL_2(\C)$ considered up to conjugation by elements in $\SL_2(\C)$, see for instance \cite{Goldman4}. 

In the above-mentioned works it has become apparent that the dimension of $\Defo(M)$ is much larger once vertices of higher valence are present, more precisely one has $\dim_\R \Defo(M) = 2N + \sum_{j=1}^k 2(m_j-3)$, where $k$ is the number of singular vertices and $m_j$ the valence of the vertex $v_j \in \Sigma$. This corresponds to the fact that the spherical cone-surface structure on the link of such a vertex becomes more flexible. More precisely, as a consequence of \cite{LuoTian,Tro}, see also \cite{MaWe}, one has that for $m \geq 3$ the space of spherical cone-manifold structures on $(\mathbb{S}^2,\{p_1, \ldots,p_m\})$ is locally parametrized by $\mathcal{T}_{0,m} \times (0,2\pi)^m$, where $\mathcal{T}_{0,m}$ is the Teichm\"uller space of the $m$-times punctured sphere. Note that the dimension of $\mathcal{T}_{0,m}$ is precisely $2(m-3)$. 

\medskip

The aim of this article is to describe how this additional flexibility can be used to deform a given cone-manifold structure on $X$ into cone-manifold structures on $X$ of possibly different topological type, i.e.~with possibly the pair $(X,\Sigma')$ not being homeomorphic to the pair $(X,\Sigma)$. This will be achieved by splitting a vertex $v_j$ of valence $m_j\geq 4$ into two or more vertices of lower valence. Note that what we actually deform is the incomplete hyperbolic metric on the smooth part $M$; the new cone-manifold is its completion.
Before describing our main result, we need to set up some more notations.

For simplicity we will assume that $\Sigma$ does not contain circle components. This makes for a cleaner statement of our main results and, besides that, the case of circle components is already well understood by the work of Hodgson and Kerckhoff. For $\varepsilon>0$ sufficiently small the subset $\bar M_\varepsilon = M \setminus U_\varepsilon(\Sigma)$ is a {\em compact core} of $M$, i.e.~$\bar M_\varepsilon \hookrightarrow M$ is a homotopy equivalence. Its boundary $\partial \bar M_\varepsilon$ is a (possibly disconnected) surface of genus  $g$, such that $\sum_{j=1}^k (m_j-3) + N = 3g-3$ (recall that $N$ and $k$ denote respectively the number of edges and of vertices contained in the singular locus $\Sigma$, and that for each vertex $v_j$, $m_j$ denotes the number of edges meeting in $v_j$). Let $L_j$ be the link of the vertex $v_j$; it is a spherical cone-surface, homeomorphic to the 2-sphere when the cone angles are less than $2\pi$. Its smooth part is an $m_j$-times punctured sphere, denoted by $N_j$, which embeds naturally (but not isometrically) as an open subsurface of $\partial \bar M_\varepsilon$; throughout this article $N_j$ and the image of this embedding will be identified. 

\begin{figure}
\centering 
\includegraphics[scale=0.9]{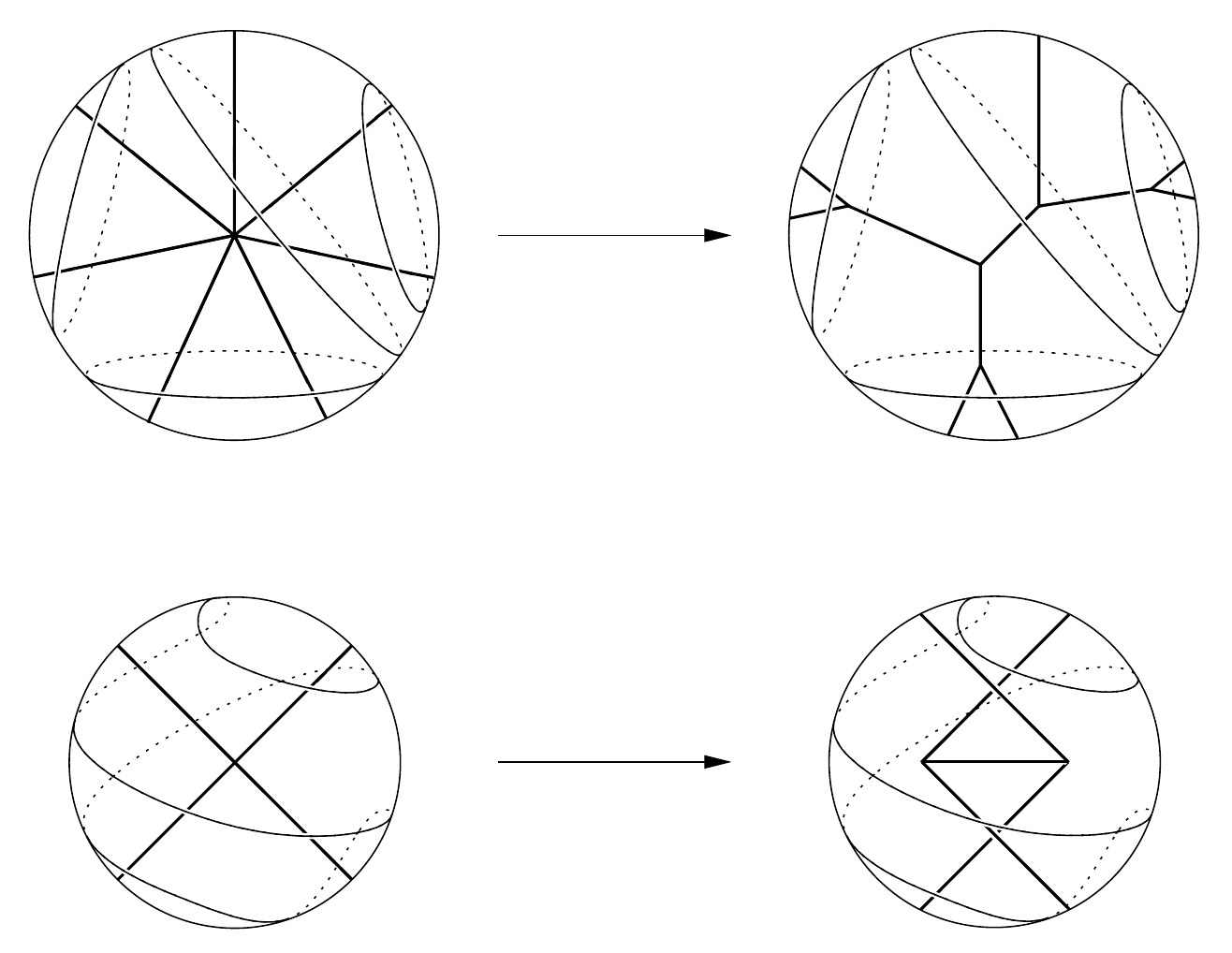}
\caption{Different ways of splitting a vertex}
\label{fig:exsplit}
\end{figure}

A pair-of-pants decomposition $\mathcal{C}_j=\{\nu_{j,1}, \ldots, \nu_{j,m_j-3}\}$ of $N_j$ determines an unknotted embedding of a trivalent tree $\Sigma'_j$ into the closed 3-ball $\overline{B_\varepsilon(v_j)}$ in such a way that the curves $\nu_{j,i}$ are precisely the meridian curves of the newly created edges $e'_{j,i}$ and that $\Sigma \cap \overline{B_\varepsilon(v_j)}$ is obtained back by collapsing these new edges. (An embedding of a trivalent tree $(T,\partial T) \hookrightarrow (D^3,\partial D^3)$ is {\em unknotted}, if it factors through an embedding of the 2-disk $(D^2,\partial D^2) \hookrightarrow (D^3,\partial D^3)$.) If we now replace $\Sigma \cap \overline{{B}_\varepsilon(v_j)}$ by $\Sigma'_j$, then we say that $(X,\Sigma')$ is obtained from $(X,\Sigma)$ by {\em splitting a vertex}, see Fig.~\ref{fig:exsplit} for examples. The homeomorphism type of the pair $(X,\Sigma')$ is determined by the pair-of-pants decomposition $\mathcal C_j$, where homotopic pair-of-pants decompositions clearly yield the same type. 

In order to construct a hyperbolic cone-manifold structure on the pair $(X,\Sigma')$, even locally on a neighbourhood of the new singular locus $\Sigma'$, we need an additional, geometric condition on the splitting curve $\nu_{j,i}$, namely that it satisfies the so-called {\em splitting condition}, cf.~Definition \ref{def:splitting_condition}. We will also simply say that such a curve is {\em splittable}. The splitting condition allows us to construct a concrete model for the splitting deformation along $\nu_{j,i}$; this construction is carried out in section \ref{sec:split}. The geometry of a splitting deformation will in general be different for two splittable curves in the same homotopy class of a simple closed curve. However, there is a notion of equivalence of splittable curves in a homotopy class, which is essentially requiring the curves to be homotopic through splittable curves, cf.~section \ref{sec:splitcond}. Equivalent splittable curves then turn out to yield the same geometric splitting deformations.

Let now $\vec\mu = \{\mu_1,\ldots, \mu_N\}$ be the set of meridians of
$\Sigma$, and let further $\mathcal{C}_j$
be a pair-of-pants decomposition of $N_j$ for each $j=1, \ldots,k$. Note that
we fix a system of simple closed curves here and not just their homotopy
classes. Let $\vec \nu = \{\nu_1, \ldots, \nu_{3g-3-N}\}$ denote the family of curves $\bigcup_{j=1}^k \mathcal{C}_j$, which is a
pair-of-pants decomposition of $\coprod_{j=1}^k N_j$; and let $\mathcal{C}$ be equal to
$\vec\mu \cup \vec\nu$, which is (up to homotopy) a pair-of-pants decomposition of $\partial
\bar M_\varepsilon$. Let us assume for simplicity now that all the curves in
$\vec \nu$ satisfy the splitting condition. We will then say that a
deformation of the hyperbolic cone-manifold structure on $X$ is {\em
  compatible} with $\vec\nu$, if its singular locus $\Sigma'$ is obtained by
splitting some vertices $v_j$ as determined by $\mathcal{C}_j$ and the
geometric deformation on a neighbourhood of the singular locus is precisely
given by a model deformation as in Definition \ref{def:nucomp}. Here
we allow that some of the newly created edges have length $0$. These
model deformations combine the splitting deformations along the curves in
$\vec \nu$ with the usual deformations of the singular tube, namely changing
the spherical structure on the links (which involves changing the cone-angles)
and changing the length or twist parameters of the singular edges. Let $C_{-1}(X,\vec\nu)$ denote the space of $\vec\nu$-compatible hyperbolic cone-manifold structures on $X$, cf.~Definition \ref{def:compatible}; it is topologized as a subset of $\Defo(M)$. Clearly $C_{-1}(X,\Sigma) \subset C_{-1}(X,\vec\nu)$ for any such family $\vec\nu$. 

To a cone-manifold structure on $X$ compatible with $\vec\nu$, we can associate the vector
of the cone angles $\alpha=(\alpha_1, \ldots, \alpha_N)$ corresponding to the
original edges $(e_1, \ldots,e_N)$ and the vector of the lengths
$\ell=(\ell_1, \ldots, \ell_{3g-3-N})=(l_{j,i})_{j=1, \ldots, k,\,i=1, \ldots,m_j}$ of the newly created edges. Our main result is that these data actually provide a parametrization of the space $C_{-1}(X,\vec\nu)$ near $X$, under the following two assumptions. The first one is that the curves in $\vec\nu$ satisfy the {\em splitting condition} as introduced above. The second one is that $\mathcal{C}$ is {\em admissible} in the sense of Definition \ref{def:admissibility}; this is essentially an algebraic assumption which ensures that the holonomy representation restricted to any pair-of-pants in a decomposition of a link remains irreducible. 

\begin{thm}\label{thm:main} 
Let $X$ be a hyperbolic cone-$3$-manifold with cone angles less than $2\pi$ and
meridian set $\vec \mu$. Let $\vec\nu$ be a pair-of-pants decomposition of $\coprod_{j=1}^kN_j$ such that ${\mathcal C}= \vec \mu \cup \vec \nu$ gives an admissible pair-of-pants decomposition of $\partial \bar M_\varepsilon$. If all the curves in $\vec\nu$ are splittable, then the map
\begin{equation*}
(\alpha,\ell): C_{-1}(X, \vec\nu) \rightarrow (0,2\pi)^N \times \R_{\geq 0}^{3g-3-N}
\end{equation*}
sending a $\vec\nu$-compatible cone-manifold structure to the vector composed of its original edges' cone angles and new edges' lengths, is a local homeomorphism at the given structure.
\end{thm}

\noindent Note that we recover Theorem \ref{thm:former} by setting $\ell_i =0$ for all $i$ (if there exists a family $\vec{\nu}$ satisfying the above assumptions).
\medskip

The main difficulty in the proof of Theorem \ref{thm:main} is to obtain an adapted local chart on $\Defo(M)$, or equivalently (using \eqref{eq:hol}) on $X(\pi_1M,\SL_2(\C))$. 
As in the previous works, the strategy consists of constructing first a larger, adequate coordinate system on the simpler space $X(\pi_1\partial \bar M_\varepsilon,\SL_2(\C))$. 
Using the infinitesimal rigidity, we can then show that part of these coordinates lift via the natural map $X(\pi_1M,\SL_2(\C)) \to X(\pi_1\partial \bar M_\varepsilon,\SL_2(\C))$ to a local chart on the former space. This is what is done in section \ref{sec:coord}. 
The coordinate system is provided by action-angle variables: as explained in subsection \ref{sec:complex-sympl} and \ref{sec:actionangle}, the character variety $X(\pi_1\partial \bar M_\varepsilon,\SL_2(\C))$ has a canonical complex-symplectic structure (actually introduced by Goldman \cite{Goldman5}, along the lines of \cite{AtBo,Goldman1}), and the traces of the curves in $\mathcal{C}$ yield a holomorphic completely integrable system whenever $\mathcal C$ is admissible. 
The $L^2$-cohomology vanishing result of \cite{We2} is then applied in subsection \ref{sec:L2} to show which part of the action-angle coordinates lifts well to $\Defo(M)$ near $X$ (Theorem \ref{thm:param}).

To construct actual deformations of the cone-manifold $X$, we follow the same
strategy of beginning with a simpler problem, namely deforming a neighbourhood
of the singular locus. The splitting condition, explained in subsection
\ref{sec:splitcond}, ensures the existence of splitting deformations as
constructed in subsection \ref{sec:split}. We survey the other, more standard deformations of the singular tube $U_\varepsilon(\Sigma)$ in subsection \ref{sec:defandhol}, and give in Proposition \ref{prop:defoU} the relation between the $\vec{\nu}$-compatible deformations of $U_\varepsilon(\Sigma)$ and the action-angle coordinates on $X(\pi_1\partial \bar M_\varepsilon,\SL_2(\C))$. 

Combining the results of sections \ref{sec:coord}  and \ref{sec:defosingtube}, we can characterize the holonomy representations of the elements of $C_{-1}(X, \vec\nu)$. This enables us to determine which $\vec\nu$-compatible deformations of $U_\varepsilon(\Sigma)$ can be extended to the whole of the cone-manifold, leading to the following statement:

\medskip

\noindent {\bf Theorem \ref{thm:construction}.}
{\em Let $X$ be a hyperbolic cone-$3$-manifold with meridian set $\vec{\mu}$ and cone angles $\alpha = (\alpha_1,\ldots,\alpha_N) \in (0,2\pi)^N$. Let $\vec\nu$ be a family of curves on $\coprod_{j=1}^kN_j$ such that ${\mathcal C}= \vec \mu \cup \vec \nu$ is up to homotopy an admissible pair-of-pants decomposition of $\partial \bar M_\varepsilon$. Then the $\vec{\nu}$-compatible deformation of the hyperbolic cone-manifold structure on $X$ corresponding to a vector $(\alpha',\ell) \in (0,2\pi)^N \times \R_{\geq 0}^{3g-3-N}$ close to $(\alpha,0)$ exists, if the curves $\nu_i$ with $\ell_i>0$ are splittable and disjoint.}

\medskip

Theorem \ref{thm:main} is then a direct consequence. In the remaining part of  section \ref{sec:localshape}, we discuss the stratified structure of the space of cone-manifold deformations of $X$ and its relation with the curve complex of $\coprod N_j$, and we explain how our main result applies to the polyhedral case.

\bigskip

\noindent {\bf Acknowledgments:} The second author would like to thank Bill Goldman for helpful conversations about surface group representations and the hospitality during the visit at University of Maryland.

\section{Local coordinates on the variety of characters} \label{sec:coord}

We begin by recalling various facts about the variety of representations of a surface group into the Lie group $\SL_2(\C)$. Let $S$ be a closed orientable surface of genus $g \geq 2$ in the following. Let $R(\pi_1S,\SL_2(\C))$ be the space of group homomorphisms $\rho: \pi_1S \rightarrow \SL_2(\C)$ equipped with the compact-open topology. The group $\SL_2(\C)$ acts on the space $R(\pi_1S,\SL_2(\C))$ by conjugation and we may form the set-theoretic quotient
$$
X(\pi_1S,\SL_2(\C))=R(\pi_1S,\SL_2(\C)) / \SL_2(\C).
$$
We will endow $X(\pi_1S,\SL_2(\C))$ with the quotient topology and we will refer to the elements of $X(\pi_1S,\SL_2(\C))$ as characters. In our main application $S$ will be the boundary of a compact core of the smooth part of a hyperbolic cone-$3$-manifold, i.e.~$S=\partial \bar M_\varepsilon$. Note that here $S$ may be disconnected, in which case we have to replace $X(\pi_1S,\SL_2(\C))$ by the product of the character varieties of the fundamental groups of the connected components of $S$, as for example explained in \cite{Montcouq2} and \cite{Weg}. Since this presents no further technical difficulty, we will mostly pretend that $S$ is actually connected.

Recall that a representation $\rho: \pi_1S \rightarrow \SL_2(\C)$ is called irreducible if it does not leave invariant any line in $\C^2$. Let $R_{irr}(\pi_1S,\SL_2(\C))$ be the space of irreducible representations and $X_{irr}(\pi_1S,\SL_2(\C))$ the space of irreducible characters. It is known that $R_{irr}(\pi_1S,\SL_2(\C))$ is a complex manifold of complex dimension $6g-3$, cf.~\cite{Goldman1,Goldman5}. Furthermore $\PSL_2(\C)$ acts properly and freely on $R_{irr}(\pi_1S,\SL_2(\C))$, hence $X_{irr}(\pi_1S,\SL_2(\C))$ is a complex manifold of complex dimension $6g-6$, see \cite{Goldman5} and the references therein for details.

\subsection{The complex symplectic structure and complex twist flows}\label{sec:complex-sympl}

In the following we review Goldman's construction of a natural complex-symplectic structure on $X_{irr}(\pi_1S,\SL_2(\C))$ as carried out in \cite{Goldman1,Goldman5}. Recall that a complex-symplectic manifold $Y$ is a complex manifold equipped with a non-degenerate closed holomorphic 2-form $\Omega$. If $2n$ is the complex dimension of $Y$, then non-degeneracy of $\Omega$ may be rephrased by saying that $\Omega^n$ is a non-vanishing $(2n,0)$-form. Let
$$
b: \sl_2(\C) \times \sl_2(\C) \rightarrow \C
$$
be a non-degenerate symmetric bilinear form which is $\Ad$-invariant. Such a pairing is provided by the trace form on $\sl_2(\C)$, i.e.~$b(A,B)=\tr (AB)$ for $A,B \in \sl_2(\C)$. Let
$$
\E_\rho = \tilde S \times_{\Ad \circ \rho } \sl_2(\C)
$$
be the flat vector bundle on $S$ with fiber $\sl_2(\C)$ and holonomy given by the representation $\Ad \circ \rho$. Then, by $\Ad$-invariance, $b$ induces a pairing
$$
b^{\E_\rho}: \E_\rho \times \E_\rho \rightarrow \C
$$
which is parallel and fiberwise non-degenerate. Finally, by Poincar\'e-duality, the pairing
\begin{align*}
\Omega : H^1(S;\E_\rho) \times H^1(S;\E_\rho) &\rightarrow \C\\
([\alpha],[\beta]) & \mapsto \int_S b^{\E_\rho}(\alpha \wedge \beta)
\end{align*}
is skew-symmetric and non-degenerate. Now recall that, using Weil's construction, the tangent space of $X_{irr}(\pi_1S,\SL_2(\C))$ at a character $[\rho]$ may be identified with the cohomology group $H^1(\pi_1S; \sl_2(\C)_{\Ad \circ \rho})$, which in turn may be identified with $H^1(S;\E_\rho)$ using de Rham's theorem. With these identifications, $\Omega$ becomes a $2$-form, shown to be indeed closed and holomorphic in \cite{Goldman1}. Hence it defines a complex-symplectic structure on $X_{irr}(\pi_1S,\SL_2(\C))$, which is actually $\pi_1S$-invariant.

\medskip

We continue with the description of the complex Hamiltonian flows associated to the trace functions; the material we present here is in essence contained in \cite{Goldman2}, see also \cite{Goldman5}. If $(Y, \Omega)$ is a complex-symplectic manifold and $f: Y \rightarrow \C$ is holomorphic, then the complex Hamiltonian vector field $X_f$ is defined by
$$
\Omega(X_f, \,\cdot\,) = df.
$$
By definition, $X_f$ is a holomorphic vector field on $Y$, and hence possesses a local holomorphic flow $\phi^f$, the {\em complex Hamiltonian flow} associated with $f$. The trajectory through $y \in Y$ is the holomorphic curve $z \mapsto \phi^f_z(y)$ satisfying the ODE
$$
\frac{\partial}{\partial z} \phi_z^f(y) = X_f(\phi_z^f(y))
$$
in ``complex time'' $z \in \C$ (for $|z|$ sufficiently small). For holomorphic functions $f,g : Y \rightarrow \C$ the Poisson bracket is as usual defined by
$$
\{f,g\} = \Omega(X_f,X_g)
$$
and we say that $f$ and $g$ Poisson-commute if $\{f,g\}=0$.

Let now $Y=X_{irr}(\pi_1S,\SL_2(\C))$. For any closed curve $\gamma$ on $S$, we define the function
\begin{align*}
\tr_\gamma: X(\pi_1 S, \SL_2(\C)) &\rightarrow \C\\
\chi = [\rho] & \mapsto \tr \rho(\gamma_0)
\end{align*}
where $\gamma_0 \in \pi_1 S$ is freely homotopic to $\gamma$; note that it does not depend on its orientation.
In the case $\gamma$ is a {\em simple} closed curve, the associated complex Hamiltonian vector field $X_\gamma$ can be described as follows.

For $A \in \SL_2(\C)$ we may consider the differential of the trace at $A$ as a map $(d \tr)_A: \sl_2(\C) \rightarrow \C$ by setting
$$
(d \tr)_A(B) = \textstyle{\left.\frac{d}{dt} \right\vert_{t=0}} \tr (A \exp(tB))
$$
for $B \in \sl_2(\C)$. If $b:\sl_2(\C) \times \sl_2(\C) \rightarrow \C$ is a non-degenerate bilinear form, the {\em variation function} $F: \SL_2(\C) \rightarrow \sl_2(\C)$ with respect to $b$ is defined by requiring that
$$
(d \tr)_A(B) = b(F(A),B)
$$
for all $B \in \sl_2(\C)$. If $b$ is chosen to be the trace form, i.e.~$b(A,B)=\tr(AB)$ for $A,B \in \sl_2(\C)$, one obtains that $F(A) = A - \frac{1}{2}\tr A \cdot \id$ for $A \in \SL_2(\C)$. Note that by virtue of this formula $\Ad(A)F(A)=F(A)$, i.e.~$F(A) \in \sl_2(\C)$ is $\Ad(A)$-invariant. 

If we identify the fiber of $\E_\rho$ over $\gamma(0)$ with $\sl_2(\C)$, then the $\Ad(\rho(\gamma))$-invariant element $F(\rho(\gamma)) \in \sl_2(\C)$ defines a parallel section $\sigma_\gamma$ on a collar $C = [0,1] \times \mathbb{S}^1$ to the left of $\gamma$, i.e.~$\partial C = \gamma'^{-1} \cup \gamma$, taking into account the orientations of $\gamma$ and the surface $S$. Let $\varphi: [0,1] \rightarrow [0,1], z \mapsto \varphi(z)$ be a smooth function which vanishes near $z=0$ and is identically $1$ near $z=1$. We set
$$
\omega_\gamma= d\varphi \otimes \sigma_\gamma \in \Omega^1(S;\E_\rho).
$$
Let $X_\gamma=[\omega_\gamma]\in H^1(S;\E_\rho)$. We claim that $\Omega(X_\gamma,[\eta]) = d\tr_\gamma([\eta])$ for all closed $\eta \in \Omega^1(S,\E_\rho)$. Note that $b^{\E_\rho}(\omega \wedge \eta) = d (\varphi\cdot b^{\E_\rho}(\sigma_\gamma \wedge \eta))$ since $\sigma_\gamma$ and $b^{\E_\rho}$ are parallel. Applying Stokes' theorem we obtain
$$
\int_S b^{\E_\rho}(\omega_\gamma \wedge \eta) = \int_{\partial C} \varphi \cdot b^{\E_\rho}(\sigma_\gamma \wedge \eta) = \int_\gamma b^{\E_\rho}(\sigma_\gamma \wedge \eta).
$$
On the other hand
$$
d\tr_\gamma([\eta])=(d \tr)_{\rho(\gamma)}( {\textstyle{\int_\gamma}} \eta) = b(F(\rho(\gamma)), {\textstyle{\int_\gamma}} \eta) =\int_\gamma b^{\E_\rho}(\sigma_\gamma \wedge \eta)
$$
where $\int_\gamma \eta \in \sl_2(\C)$ is defined using parallel transport in $\E_\rho$ along $\gamma$, cf.~\cite{We1} for details. The claim follows.

Let now $\mathcal{C}=\{\gamma_1, \ldots, \gamma_{3g-3}\}$ be a pair-of-pants decomposition of $S$. We may arrange that the supports of the forms $\omega_{\gamma_1}, \ldots,\omega_{\gamma_{3g-3}}$ are disjoint, hence we immediately deduce the following result of Goldman, cf.~Prop.~2.2.2 and Cor.~2.2.3 in \cite{Goldman5}:

\begin{prop}\label{prop:poisson-commute}
The functions $\tr_{\gamma_1}, \ldots, \tr_{\gamma_{3g-3}}$ Poisson-commute pairwise, i.e.~$\Omega(X_{\gamma_i},X_{\gamma_j})=0$ for all $1 \leq i,j \leq 3g-3$.  
\end{prop}

In particular, the complex Hamiltonian flows $\phi_z^i$ associated with the functions $\tr_{\gamma_i}$ commute and define a local holomorphic $\C^{3g-3}$-action $$\phi_{z} = \phi^1_{z_1} \circ \ldots \circ \phi^{3g-3}_{z_{3g-3}}
$$
on $X_{irr}(\pi_1S,\SL_2(\C))$, where $z = (z_1, \ldots, z_{3g-3}) \in \C^{3g-3}$ (for $\|z\|$ sufficiently small).
Under the assumption that the differentials $d\tr_{\gamma_1}, \ldots, d\tr_{\gamma_{3g-3}}$ are linearly independent at $\chi_0 \in  X_{irr}(\pi_1S,\SL_2(\C))$, one thus obtains a holomorphic completely integrable system
$$
\tr_{\vec \gamma} = (\tr_{\gamma_1}, \ldots, \tr_{\gamma_{3g-3}}) : X_{irr}(\pi_1S,\SL_2(\C))  \supset \mathcal U \rightarrow \C^{3g-3}.
$$
By choosing a local section
$$
\sigma : \C^{3g-3} \supset \mathcal V \rightarrow \tr_{\vec \gamma}^{-1}(\mathcal V) \subset X_{irr}(\pi_1S,\SL_2(\C)) 
$$
with $\sigma(\tr_{\vec \gamma}(\chi_0))=\chi_0$ we obtain local ``action-angle'' coordinates
\begin{align*}
 X_{irr}(\pi_1S,\SL_2(\C)) \supset \mathcal U & \rightarrow \C^{3g-3} \times \C^{3g-3}\\
\chi & \mapsto (\tr_{\vec \gamma}(\chi), \tau_{\vec \gamma}(\chi))
\end{align*}
near $\chi_0$; they are defined by requiring that $\chi=\phi_{\tau_{\vec \gamma}(\chi)}(\sigma(\tr_{\vec \gamma}(\chi))$ for $\chi \in \mathcal{U}$ and using the implicit function theorem. This set of coordinates will be further discussed in section \ref{sec:actionangle}.

As explained by Goldman in \cite{Goldman5}, the complex Hamiltonian flow $\phi_z$ on $X_{irr}(\pi_1S, \SL_2(\C))$ associated with the function $\tr_{\gamma}$ is covered by a holomorphic flow $\psi_z$ on $R_{irr}(\pi_1S,\SL_2(\C))$. This {\em complex twist flow} can be quite explicitly described as follows (we have to be rather careful about various choices in order to match the directions of the flows): 
\begin{enumerate}
\item If $\gamma$ is separating, then $\pi_1S$ splits as an amalgamated free product. Let $S_1$ denote the component of $S$ cut along $\gamma$ lying to the left of $\gamma$, taking into account the orientations of $\gamma$ and the surface $S$. Let $S_2$ denote the component lying to the right of $\gamma$. Then $\pi_1S = \pi_1 S_1 \star_{\gamma} \pi_1S_2$. 
For $\rho \in R(\pi_1S, \SL_2(\C))$ and $z \in \C$ we set
$$
\psi_z(\rho)(\gamma') = 
\left\{ 
\begin{array}{c@{\quad:\quad}c}
\rho(\gamma') & \gamma' \in \pi_1S_1\\ 
\zeta_z \rho(\gamma')\zeta_{-z} & \gamma' \in \pi_1S_2
\end{array}
\right., 
$$
where $\zeta_z$ denotes the complex 1-parameter subgroup in $\SL_2(\C)$ associated with $\rho(\gamma)$. More precisely, if $F:\SL_2(\C) \rightarrow \sl_2(\C)$ is the variation function with respect to the trace form $b$, then $\zeta_z = \exp (z F(\rho(\gamma)))$.
Note in particular that $\zeta_z$ centralizes $\rho(\gamma)$ in $\SL_2(\C)$. To check that the flows match, we calculate the periods of $\omega_\gamma$:
\begin{align*}
\int_{\gamma'} \omega_\gamma &=
\left\{ 
\begin{array}{c@{\quad:\quad}c}
0 & \gamma' \in \pi_1S_1\\ 
(1-(\Ad \circ \rho)(\gamma')) F(\rho(\gamma)) & \gamma' \in \pi_1S_2
\end{array}
\right.\\
&=\textstyle{\left.\frac{d}{dz} \right\vert_{z=0}}\psi_z(\rho)(\gamma')\rho(\gamma')^{-1}. 
\end{align*}
\item If $\gamma$ is non-separating, then $\pi_1S$ splits as an HNN-extension. Let $S'$ denote the connected surface obtained by cutting $S$ along $\gamma$. Then $\pi_1S = \pi_1 S' \star_{\gamma}$.
More precisely, let $\lambda$ be another simple closed curve intersecting $\gamma$ transversally at the base-point with positive intersection number. This determines an arc, again denoted by $\lambda$, in $S'$ connecting $\partial_+S'$ to $\partial_-S'$. If we place the base-point on $\partial_-S'$ and denote the element in $\pi_1S'$ corresponding to $\gamma$ by $\gamma_-$, then with $\gamma_+=\lambda^{-1} \star \gamma \star \lambda$ we get $\pi_1 S = \langle \pi_1 S', \lambda \vert \lambda^{-1} \gamma_- \lambda= \gamma_+ \rangle$. For $\rho \in R(\pi_1S, \SL_2(\C))$ and $z \in \C$ we set 
$$
\psi_z(\rho)(\gamma') = 
\left\{ 
\begin{array}{c@{\quad:\quad}c}
\rho(\gamma') & \gamma' \in \pi_1S'\\ 
\zeta_{z}\rho(\lambda) & \gamma'=\lambda
\end{array}
\right., 
$$
with $\zeta_z = \exp (z F(\rho(\gamma)))$ as above. Again as a check we calculate the periods of $\omega_\gamma$:
\begin{align*}
\int_{\gamma'} \omega_\gamma &=
\left\{ 
\begin{array}{c@{\quad:\quad}c}
0 & \gamma' \in \pi_1S'\\ 
F(\rho(\gamma)) & \gamma' = \lambda
\end{array}
\right.\\
&=\textstyle{\left.\frac{d}{dz} \right\vert_{z=0}}\psi_z(\rho)(\gamma')\rho(\gamma')^{-1}. 
\end{align*}
\end{enumerate}
The significance of this for us is that we can easily reconstruct $\chi \in \mathcal{U}$ from its coordinates $\tr_{\vec \gamma}(\chi)$ and $\tau_{\vec \gamma}(\chi)$ once we know the local section $\sigma$.

\subsection{Action-angle coordinates} \label{sec:actionangle}

We return now to our $3$-dimensional situation, i.e.~$X$ is a closed, orientable hyperbolic cone-$3$-manifold with smooth part $M = X\setminus \Sigma$ and $S=\partial \bar M_{\varepsilon}$; we will furthermore assume that all the cone angles are smaller than $2\pi$. The hyperbolic metric on $M$ determines a holonomy representation $\hol : \pi_1M \to \SL_2(\C)$, which in turns induces on $\partial \bar M_{\varepsilon}$ a representation $\rho_0 \in R(\pi_1 \partial \bar M_{\varepsilon}, \SL_2(\C))$, with corresponding character $\chi_0=[\rho_0]$. 

For each edge $e_i$ of the singular locus, let $\mu_i \subset \partial \bar M_{\varepsilon}$ be a meridian of $e_i$, i.e.~a simple closed curve winding exactly once around $e_i$. The collection $\vec \mu = \{\mu_1,\ldots, \mu_N\}$ decomposes $\partial \bar M_{\varepsilon}$ into a family of $k$ subsurfaces $\{\Sigma_1,\ldots, \Sigma_k\}$, each of which is homeomorphic to the regular part $N_j$ of the link of a vertex $v_j$; in the following we will always identify $\Sigma_j$ and $N_j$.

As in \cite{Montcouq2,We2} we find the following (see \cite{Montcouq2} for a detailed proof, where this result appears as Theorem 6):

\begin{lemma}
  $\hol: \pi_1 N_j \rightarrow \SL_2(\C)$ is irreducible.
\end{lemma}

Actually, the induced holonomy on $N_j$ fixes a point in $\H^3$ (corresponding to the vertex $v_j$), so has values in a maximal compact subgroup $K$ of $\SL_2(\C)$; up to conjugacy we can assume $K = \SU(2)$.

For each $1\leq j \leq k$, let $\mathcal{C}_j$ be a pair-of-pants decomposition of $N_j$. If we cut $N_j$ along $\mathcal{C}_j$ we obtain a disjoint collection $P_{j,1}, \ldots, P_{j,m_j-3}$ of subsurfaces, each homeomorphic to a thrice punctured sphere.

\begin{Def}\label{def:admissibility} Let $S$ be a surface and $\rho : \pi_1S \to \SL_2(\C)$ an irreducible representation. 
A pair-of-pants decomposition $\mathcal{C}$ of $S$ is called $\rho$-{\em admissible} (or just admissible for short) if the restriction of $\rho$ to each of the pairs of pants obtained by cutting $S$ along $\mathcal{C}$ is irreducible.
\end{Def}

A natural question is whether an admissible pair-of-pants decomposition for $N_j$ always exists. This is settled by the following:

\begin{prop}\label{prop:admdecomp}
Let $S$ be a $d$-punctured sphere and $\rho : \pi_1S \to \SU(2)$ an irreducible representation, such that the holonomy of any peripheral element is non-trivial. Then $S$
admits a $\rho$-admissible pair-of pants decomposition.
\end{prop}

\begin{pf}
The proof is by induction on $d$. If $d=3$, then the result is trivial, so we will now assume that $d>3$.
Since $\rho$ has values in $\SU(2)$, it is reducible if and only if all the elements of its image commute. Let
$\langle \gamma_1,\ldots , \gamma_d\ |\ \prod_i \gamma_i =1 \rangle$ be a presentation of $\pi_1S$ such that
$\gamma_i$ is a loop around the $i$-th puncture of $S$. Since $\rho$ is irreducible, up to a change in the order of
the generators we can assume that $\rho(\gamma_1)$ and $\rho(\gamma_2)$ do not commute. The curve $\gamma_1.\gamma_2$ (or rather, any simple closed curve homotopic to $\gamma_1.\gamma_2$) 
cuts $S$ into a pair-of-pants $P'$ and a $(d-1)$-punctured sphere $S'$. The induced representation on $P$ is
generated by $\rho(\gamma_1)$ and $\rho(\gamma_2)$ and is thus irreducible. If the induced representation on $S'$ is
also irreducible, then by induction we can find an admissible decomposition of $S'$ and the proof is completed. If on
the other hand the induced representation on $S'$ is reducible, this means that there exists a complex line $L$ which
is invariant by $\rho(\gamma_i)$ for all $3\leq i \leq d$. If $L$ is invariant by $\rho(\gamma_2)$, then it is also
invariant by $\rho(\gamma_1) = \left(\rho(\gamma_2) \ldots \rho(\gamma_d) \right)^{-1}$, and this contradicts the
irreducibility of $\rho$. So $L$ is not invariant by $\rho(\gamma_2)$; similarly, it is not invariant by
$\rho(\gamma_1)$. Now any simple closed curve homotopic to $\gamma_2.\gamma_3$ cuts $S$ into a pair-of-pants $P''$ and a $(d-1)$-punctured sphere
$S''$. The induced representation on $P''$ is generated by $\rho(\gamma_2)$ and $\rho(\gamma_3)$, and they do not
commute since $L$ is invariant by $\rho(\gamma_3)$, which is not trivial, but not by $\rho(\gamma_2)$; this implies that the induced representation on
this pair-of-pants is irreducible. Similarly, the induced representation on $S''$ is irreducible since it is
generated by $\rho(\gamma_1), \rho(\gamma_4),\ldots, \rho(\gamma_d)$, which do
not commute.\end{pf}

The same argument shows that under the same assumptions, any simple closed curve $\gamma \subset S$ with the property that
$\rho(\gamma)$ is non-trivial and that the restriction of $\rho$ to each component
obtained by cutting $S$ along $\gamma$ is irreducible, can be completed to an
admissible pair-of-pants decomposition of $S$. In the case where $S=N_j$ and $\rho=\rho_0$, the peripheral elements are actually edges' meridians whose holonomy is non-trivial since the cone angles are smaller than $2\pi$, so the proposition applies.

\begin{lemma}\label{submersion_link}
Let $S$ be a $d$-punctured sphere and $\rho : \pi_1S \to \SL_2(\C)$ an irreducible representation. Let $\{\gamma_1,\ldots,\gamma_d\}$ be the boundary curves of $S$ and assume that the pair-of-pants decomposition $\mathcal{C}=\{\gamma_{d+1}, \ldots, \gamma_{2d-3}\}$ is admissible. Then the map
$$
(\tr_{\gamma_1}, \ldots, \tr_{\gamma_{2d-3}}): X_{irr}(\pi_1 S, \SL_2(\C)) \rightarrow \C^{2d-3}
$$
is a submersion at $\rho$.
\end{lemma}

The proof of this lemma is based on the following elementary fact (see for instance \cite{Goldman3}, p.~578):

\begin{lemma}\label{submersion_pop}
Consider the map 
\begin{align*}
f: \SL_2(\C) \times \SL_2(\C) &\rightarrow \C^3\\
(A,B) & \mapsto (\tr A, \tr B, \tr AB).
\end{align*}
Then $df_{(A,B)}$ is surjective if and only if $A$ and $B$ do not commute. 
\end{lemma}

Let now $P$ be a thrice punctured sphere and 
$\pi_1P = \langle \gamma_1, \gamma_2, \gamma_3 \,\vert\, \gamma_1 \gamma_2 \gamma_3 \rangle$ 
a presentation of its fundamental group. Lemma \ref{submersion_pop} may be translated into the following equivalent statement:
\begin{cor}\label{cor_submersion_pop}
Consider the map
$$
\tr_{\vec\gamma}=(\tr_{\gamma_1}, \tr_{\gamma_2}, \tr_{\gamma_3}) : X(\pi_1P, \SL_2(\C)) \rightarrow \C^3.
$$
Then $(d\tr_{\vec\gamma})_\rho$ is surjective if and only if $\rho$ is nonabelian.
\end{cor}

We can then use the standard gluing construction as explained in \cite{Montcouq2} and \cite{We1,We2} together with Corollary \ref{cor_submersion_pop} to finish the proof of Lemma \ref{submersion_link}.

\medskip

Recall that for each $1\leq j \leq k$, $\mathcal{C}_j$ is a pair-of-pants decomposition of the smooth part $N_j$ of the $j$-th vertex's link, and that the set $\vec \mu =\{\mu_1,\ldots, \mu_N\}$ consists of all the edges' meridians. Let $\vec \nu = \bigcup_{j=1}^k \mathcal{C}_j = \{\nu_1,\ldots,\nu_{3g-3-N}\}$ (where $g$ stands for the genus of $\partial \bar M_\varepsilon$, so that $\sum_{j=1}^k(m_j-3) = 3g-3-N$, with $m_j$ being the valence of the $j$-th vertex). Then the collection $\mathcal{C} = \vec \mu \cup \vec \nu$ is a pair-of-pants decomposition of $\partial \bar M_\varepsilon$, which is admissible if and only if all the $\mathcal{C}_j$ are admissible. To this decomposition we associate the two trace maps
$$
\tr_{\vec \mu}=(\tr_{\mu_1}, \ldots, \tr_{\mu_N}) : X(\pi_1 \partial\bar M_\varepsilon,\SL_2(\C)) \rightarrow
\C^N
$$
and
$$
\tr_{\vec \nu}=(\tr_{\nu_1}, \ldots,\tr_{\nu_{3g-3-N}}) : X(\pi_1 \partial \bar M_\varepsilon,\SL_2(\C)) \rightarrow\C^{3g-3-N}.
$$

\begin{prop}\label{submersion}
Assume that $\mathcal{C} = \vec \mu \cup \vec \nu$ is $\rho_0$-admissible. Then the map 
$$
\Tr_\mathcal{C} = (\tr_{\vec \mu},\tr_{\vec \nu}):  X_{irr}(\pi_1 \partial\bar M_\varepsilon,\SL_2(\C)) \rightarrow
\C^{3g-3}
$$
is a submersion at $\chi_0$.
\end{prop}

The proof uses the standard gluing construction as in Lemma \ref{submersion_link}. Further details are left to the reader.

\bigskip

If $\mathcal{C} = \vec \mu \cup \vec \nu$ is admissible for $\rho_0$, then upon choosing a local section for the map $(\tr_{\vec\mu},\tr_{\vec\nu})$, we obtain using Proposition \ref{submersion} together with Proposition \ref{prop:poisson-commute} and the ensuing discussion local action-angle coordinates near $\chi_0=[\rho_0]$:
\begin{align*}
 X_{irr}(\pi_1 \partial \bar M_\varepsilon,\SL_2(\C)) \supset \mathcal{U} & \rightarrow \C^N \times \C^{3g-3-N}  \times \C^N \times \C^{3g-3-N}\\
\chi & \mapsto (\tr_{\vec\mu},\tr_{\vec\nu},\tau_{\vec\mu},\tau_{\vec\nu})
\end{align*}

\medskip

\begin{rem}{\bf (Warning)} The ``action-angle'' terminology, while customary, is somewhat misleading in our context. Indeed, it is the action variables $\tr_{\mu_i}$ and $\tr_{\nu_i}$ that actually involve the cone angles via the trace of the holonomy of elliptic isometries. We will see that the ``angle'' variables $\tau_{\mu_i}$ and $\tau_{\nu_i}$ have little to do with the cone angles but are rather related to the lengths of the edges.
\end{rem}

In our $3$-dimensional setting we can give a more geometric description of the angle-coordinates $\tau_{\vec\mu}$ and $\tau_{\vec\nu}$. But since they depend on $\sigma$, we will first look at the most natural choice for this local section. 

Let $N_j$ be the regular part of the link of a singular vertex of $X$; we know that the induced holonomy representation on $N_j$ fixes a point $p_j \in \mathbb{H}^3$ and hence has values in a maximal compact subgroup $K_j$ of $\SL_2(\C)$. In particular, the induced holonomy representation on any of the pair-of-pants of the admissible decomposition of $N_j$ has values in this maximal compact subgroup, which is conjugated to $\SU(2)$ inside $\SL_2(\C)$. A first consequence of this fact is that $\Tr_\mathcal{C}(\chi_0)$ actually belongs to $(-2,2)^{3g-3} \subset \C^{3g-3}$.

Let $t = (t_1,\ldots,t_{3g-3}) \in \R^{3g-3}$. Let $\gamma_{j_1}, \gamma_{j_2},\gamma_{j_3} \in \mathcal{C}_j$ be the boundary curves of a pair-of-pants $P$ in the decomposition of $N_j$. 
If $t$ is close enough to $\Tr_{\mathcal{C}}(\chi_0)$, then there exists a unique (up
to conjugation) representation $\rho_{P}(t) : \pi_1P \to K_j$ such that
$\tr(\rho_{P}(t)(\gamma_{j_i})) = t_{j_i}$ for $i\in \{1,2,3\}$; and we can
proceed in the same way for all the pairs-of-pants of the decomposition. It is
then possible to glue these representations together to obtain an element
$\rho_{N_j}(t) \in R_{irr}(\pi_1 N_j, K_j)$. The representations
$\rho_{N_j}(t)$ may then again be glued together to obtain $\rho(t) \in
R_{irr}(\pi_1\partial \bar{M}_\varepsilon,\SL_2(\C))$ such that
$\Tr_{\mathcal{C}}([\rho(t)]) = t$. Such an element $\rho(t)$ is not unique, but this
shows that we can choose a neighbourhood $\mathcal V \subset \C^{3g-3}$ of
$\Tr_{\mathcal{C}}(\chi_0)$ and a local section $\sigma :\mathcal V \rightarrow
X_{irr}(\pi_1 \partial \bar{M}_\varepsilon,\SL_2(\C))$ of $\Tr_{\mathcal{C}}$ such
that for any $t \in \mathcal V \cap \R^{3g-3}$, $\sigma(t)$ is the character
of a representation $\rho(t)$ with the property that the restriction of
$\rho(t)$ to $\pi_1 N_j$ has values in a maximal compact subgroup $K_j(t)$ for
all $j=1,\ldots,k$. Note that in general we cannot achieve that $K_j(t)=K_j$
for all $t\in \mathcal V \cap \R^{3g-3}$. 

We want to understand how the complex twist flows $\psi^i_z$ act on the holonomy representation $\rho_0$ of the hyperbolic cone-manifold structure (or more generally, on a representation $\rho(t)$ as constructed above). If $A \in \SL_2(\C)$ corresponds to an elliptic isometry, we may assume w.l.o.g.~that
$$
A =\begin{bmatrix}
\lambda & 0\\
0 & \bar \lambda
\end{bmatrix}
$$
for $\lambda \in \operatorname{U}(1)$. If we set $\lambda=e^{i\alpha/2}$ for $\alpha \in \R$ (i.e.~$A$ corresponds to an elliptic isometry with (oriented) rotation angle $\alpha$), then the complex 1-parameter subgroup $\zeta_z$ corresponding to $A$ is given by
$$
\zeta_z =
\begin{bmatrix}
e^{iz \sin (\alpha/2)} & 0\\
0 & e^{-iz \sin(\alpha/2)}
\end{bmatrix}.
$$
The corresponding group of isometries preserves the axis $\delta_0=\{0\} \times \R_+ \subset \C \times \R_+$ in the upper half-space model of $\mathbb{H}^3$. In particular,
$$
\zeta_t =
\begin{bmatrix}
e^{it\sin (\alpha/2)} & 0\\
0 & e^{-it\sin (\alpha/2)}
\end{bmatrix}, t \in \R
$$
is a real 1-parameter group of rotations fixing $\delta_0$ with (oriented) rotation angle given by $2t \sin(\alpha/2)$. Similarly,
$$
\zeta_{it} =
\begin{bmatrix}
e^{-t\sin (\alpha/2)} & 0\\
0 & e^{t \sin (\alpha/2)}
\end{bmatrix}, t \in \R
$$
is a real 1-parameter group of translations along $\delta_0$ with (signed) translation length given by $-2t \sin(\alpha/2)$.

This yields a rather concrete description of representations $\rho$ with character $\chi\in\mathcal U$ in terms of action-angle coordinates. The most relevant situation for us is when $\Tr_{\mathcal{C}}(\chi) \in (-2,2)^{3g-3}$ (such representations provide candidates for holonomy representations of cone-manifold structures {\em compatible} with $\vec{\nu}$ in the sense of Definition \ref{def:compatible}). In this case, as a consequence of the particular choice of the local section $\sigma$ as described above, we obtain that vanishing of the imaginary parts of the angle coordinates $\tau_{\vec\nu}$ for a character $\chi$ implies that $\rho(\pi_1 N_j)$ is contained in a maximal compact subgroup $K_j(\rho)$ for all $j=1,\ldots,k$. This corresponds geometrically to a deformation keeping the topological type $(X,\Sigma)$ fixed. Giving the angle-coordinate $\tau_{\nu_{i}}$ a non-zero imaginary part amounts to splitting along $\nu_i$ the link $N_j$ that contains this curve into components $N_j'$ and $N_j''$ and separating the fixed points of $\rho(\pi_1N_j')$ and $\rho(\pi_1N_j'')$ in $\H^3$ by some positive distance. This corresponds geometrically to a {\em splitting deformation} as described in the following sections.

\subsection{$L^2$-cohomology}\label{sec:L2}

In the beginning of section \ref{sec:complex-sympl} we have introduced the flat $\sl_2(\C)$-bundle $\E_\rho = \tilde S \times_{\Ad \circ \rho } \sl_2(\C)$ over a surface $S$, associated to a representation $\rho$. In the same way, we can define over the smooth part $M$ of the cone-manifold $X$ a bundle 
$$
\E = \tilde M \times_{\Ad \circ \hol} \sl_2(\C), 
$$
where $\hol: \pi_1M \rightarrow \SL_2(\C)$ is the holonomy representation of the hyperbolic structure on $M$. This bundle $\E$ carries a natural flat connection $\nabla^{\E}$. It is canonically identified with the bundle of infinitesimal isometries $\mathfrak{so}(TM) \oplus TM$; this identification gives a natural metric $h^{\E}$ on $\E$ (the direct sum decomposition is not preserved by the connection, so that $h^{\E}$ is not parallel with respect to $\nabla^{\E}$). Using Weil's construction and de Rham's theorem one obtains
$$
T_{\chi_0} X(\pi_1 M,\SL_2(\C)) \cong H^1(\pi_1M;\sl_2(\C)_{\Ad \circ \hol}) \cong H^1(M;\E).
$$
Note further that using the hyperbolic metric $g$ on $M$ and the bundle metric $h^{\E}$, the $L^2$-cohomology groups $H^i_{L^2}(M;\E)$ are defined. For details concerning these constructions we refer the reader to our earlier works. If $N_j$ is the smooth part of the $j$-th link, then $\hol(\pi_1N_j)$ fixes a point $p_j \in \H^3$, i.e.~is conjugated into $\SU(2)$. As a consequence one has the following splitting:
$$
\left.\E\right\vert_{N_j} = \E_j^1 \oplus \E_j^2,
$$
where the first summand corresponds to infinitesimal rotations about $p_j$ and the second one to infinitesimal translations at $p_j$. The following has been proven in \cite{We2}:

\begin{thm}\label{L2vanishing}
Let $c \in H^1_{L^2}(M;\E)$ be a class with the
property that for all vertices $v_j$ the following holds:
$$
\left.c \right\vert_{H^1_{L^2}(N_j;\E_j^1)} = 0 \quad\text{or}\quad 
\left.c \right\vert_{H^1_{L^2}(N_j;\E_j^2)} = 0.
$$
Then $c=0$.
\end{thm}

As a consequence we obtain that the map $H^1(\bar M_\varepsilon;\E) \rightarrow H^1(\partial \bar M_\varepsilon;\E)$ is injective, so that we will identify $H^1(\bar M_\varepsilon;\E)$ with its image, and that $\dim_\C H^1(\bar M_\varepsilon;\E)=\frac{1}{2} \dim_\C H^1(\partial \bar M_\varepsilon;\E)=N+\sum_{j=1}^k (m_j-3)$, cf.~Proposition 4.7 in \cite{We2}. 

The following is the main result of this section:

\begin{thm}\label{thm:param}
Let $\mathcal{C} = \vec \mu \cup \vec \nu$ be a pair-of-pants decomposition of $\partial \bar M_\varepsilon$. Assume that $\mathcal{C}$ is admissible for $\rho_0=\hol$ and $\mathcal{U}$ is a neighbourhood of $\chi_0$ in $X_{irr}(\pi_1 \partial \bar M_\varepsilon,\SL_2(\C))$, on which action-angle coordinates are defined. Then the map
\begin{align*}
\Phi_{\mathcal{C}}: \mathcal{U} \cap X(\pi_1M,\SL_2(\C)) & \rightarrow \C^N \times \R^{3g-3-N} \times \R^{3g-3-N} \\
\chi & \mapsto ( \tr_{\vec \mu}(\chi), \Im \tau_{\vec \nu}(\chi), \Im \tr_{\vec \nu}(\chi))
\end{align*}
is a local diffeomorphism at $\chi_0$.
\end{thm}

\begin{pf-lb}
Let $L := \ker (d\Phi_{\mathcal{C}})_{\chi_0} \subset H^1(\partial \bar M_\varepsilon;\E),
$
where as usual $H^1(\partial \bar M_\varepsilon;\E)$ is identified with $T_{\chi_0}X(\pi_1 \partial \bar M_\varepsilon,\SL_2(\C))$. We claim that 
$$
L \cap H^1(\bar M_\varepsilon;\E)=\{0\},
$$
where according to the remark following Theorem \ref{L2vanishing} the space $H^1(\bar M_\varepsilon;\E)$ is identified with a subspace of $H^1(\partial \bar M_\varepsilon;\E)$. 

Let now $c \in L \cap H^1(\bar M_\varepsilon;\E)$. In particular, $c \in \ker (d\tr_{\vec \mu})_{\chi_0}$, hence by Corollary 4.17 in \cite{We2}, $c \in H^1_{L^2}(M;\E)$. Note that by Corollary 4.4 in \cite{We2} the space $H^1_{L^2}(M;\E)$ may be identified with a subspace of $H^1(\bar M_\varepsilon;\E)$. Finally, since $\Im (d\tr_{\vec\nu})_{\chi_0}(c)=\Im (d\tau_{\vec\nu})_{\chi_0}(c)=0$, the restriction of $c$ to $N_j$ is tangent to a path of characters $[\rho_t]$ with $\rho_t: \pi_1N_j \rightarrow \SU(2)$ for all $j=1, \ldots,k$. This implies that
$$
\left. c \right\vert_{H^1_{L^2}(N_j;\E_j^2)}=0
$$
for all $j=1, \ldots,k$. Hence Theorem \ref{L2vanishing} applies to yield $c=0$. It follows that $\Phi_{\mathcal C}$ is an immersion at $\chi_0$. Now 
$$
\dim_\R X(\pi_1 \bar M_\varepsilon, \SL_2(\C)) = 2N + \sum_{j=1}^k 2(m_j-3)
$$ 
such that the result follows.
\end{pf-lb}

\begin{rem}
The same argument shows that $(\tr_{\vec\mu}, \Re \tr_{\vec\nu}, \Re \tau_{\vec\nu})$ also provides a system of local coordinates near $\chi_0$. However, as we will see later, the coordinate system of Theorem \ref{thm:param} is better adapted for our purposes.
\end{rem}

\section{Deformations of the singular tube}\label{sec:defosingtube}

Let $X$ be a closed hyperbolic cone-$3$-manifold. 
The topologies of the singular locus $\Sigma$ and of $\partial \bar{M}_\varepsilon$ are closely related, apart from the obvious one-to-one correspondence between the components of $\Sigma$ and those of $\partial \bar{M}_\varepsilon$: each component of $\Sigma$ is a (compact) connected graph, hence is homotopy equivalent to a bouquet of  $g$ circles, where $g$ is precisely the genus of the corresponding component of $\partial \bar{M}_\varepsilon$.

Let $X'$ be a deformation of $X$, i.e.~$X$ and $X'$ have diffeomorphic regular parts $M$ and $M'$, and the respective hyperbolic structures on $M$ and $M'$ are close. The boundary surfaces $\partial \bar{M}_\varepsilon$ and $\partial \bar{M}'_\varepsilon$  are diffeomorphic, so by the above remark the singular loci $\Sigma$ and $\Sigma'$ are homotopy equivalent. This implies that it is possible to go from one singular locus to the other by shrinking some edges to zero and splitting some vertices. But by continuity of the lengths of the edges, if $X'$ and $X$ are close enough then the only possibility is that $\Sigma'$ is obtained from $\Sigma$ by splitting some (possibly zero) vertices, that is, $X'$ has some ``new'' singular edges as compared to $X$.

In order to study the deformations of $X$, we begin with the less obstructed problem of constructing deformations of the singular tube $U_\varepsilon(\Sigma)$. We will then see (in section \ref{sec:localshape}) when they can be extended to the whole of $M$. The possible deformations of a singular tube are quite easy to classify (cf.~\cite{Maz-GM}, section 3).
\begin{itemize}
\item We can change the lengths and/or twist parameters of the existing edges.
\item We can deform the spherical structure of the vertices' links, with or without changing the cone angles. In the latter case, this also implies changing the cone angles of the existing edges.
\item Finally, we can split some vertices to create new edges.
\end{itemize}
The last case is clearly the most complicated. Topologically, such a deformation is described by specifying the meridian of the new edge, see Fig.~\ref{fig:exsplit}, but the actual geometric description is more complicated. In fact, to ensure the existence of the splitting deformation we will need an additional assumption.

\subsection{The splitting condition}\label{sec:splitcond}
 
We now introduce a geometric condition on a simple closed curve, which we call the {\em splitting condition}. This condition complements the more algebraic condition of admissibility of a pair-of-pants decomposition in the formulation of our main results.

Let $N$ be an incomplete oriented spherical surface, for instance the smooth part of a spherical cone-surface, and let $\gamma$ be a closed curve on $N$. We can choose a base-point $x$ on  $\gamma$, so that $\gamma$ defines an element of $\pi_1(N,x)$. The spherical structure on $N$ allows us to define a developing map $dev : \tilde{N} \to \mathbb{S}^2$ and associated holonomy representation $\rho : \pi_1(N,x) \to \Isom(\mathbb{S}^2)$. The image of $\gamma$ by $\rho$ is then a direct isometry of $\mathbb{S}^2$. If we assume that this isometry $f$ is not trivial, it has exactly two fixed points, denoted by $Fix(f)$. 
 
Let $r$ and $\theta$ be polar coordinates in $\mathbb{S}^2\setminus Fix(f)$. The angular variable $\theta$ is only defined modulo $2\pi$, but the corresponding $1$-form $d\theta$ and Killing vector field $\partial / \partial \theta$ are well-defined, and are invariant by $f$. 
In particular, $d\theta$ (resp.~$\partial / \partial \theta$) descends to a $1$-form (resp.~local Killing vector field) in a neighbourhood of $\gamma$. Actually, since $d\theta$ is not defined on the whole of $\mathbb{S}^2$ (it has singularities at $Fix(f)$), the corresponding $1$-form along $\gamma$ may be singular.

\begin{Def}\label{def:splitting_condition}
With the above notations, we say that a simple closed curve on $N$ with non-trivial holonomy satisfies the splitting condition (or is splittable) if it is transverse to the $1$-form  $d\theta$ (i.e.~the tangent direction to the curve at any point does not lie in $\ker d\theta$; in particular $d\theta$ has no singularity along the curve).
\end{Def}

An important fact is that this property is open, and in particular still holds if $N$ is slightly deformed:
\begin{lemma}\label{lem:open}
Let $\gamma$ be a simple closed curve on $N$. Then the set of spherical metrics on $N$ for which $\gamma$ is splittable is open.
\end{lemma} 

Indeed, small perturbations of the spherical metric induce small perturbations of $d\theta$, so that the transversality condition is preserved. 

\medskip

On the set of splittable curves on $N$, we define a new equivalence relation refining the one given by homotopy and say that two simple closed curves on $N$ satisfying the splitting condition are equivalent if they are homotopic through splittable curves. In particular, equivalent splittable curves will play exactly the same role in the remainder of the article. The following result shows that in a given homotopy class there can be at most two distinct equivalence classes of splittable curves:

\begin{prop}\label{prop:uniqueness}
Let $\gamma_1$ and $\gamma_2$ be two homotopic splittable curves on $N$. The $1$-form $d\theta$ defined in a neighbourhood of $\gamma_1$ can be uniquely extended to a neighbourhood of $\gamma_2$. If $\int_{\gamma_1} d\theta$ and $\int_{\gamma_2} d\theta$ have the same sign, then $\gamma_1$ and $\gamma_2$ are equivalent, i.e.~homotopic through splittable curves.
\end{prop}

\begin{pf}
Let $X$ denote the locally defined vector field $\sin(r)\,\partial/\partial r$ spanning the local foliation $\ker d\theta$ along $\gamma_1$ where it is nonsingular. After choosing a base-point on $\gamma_1$, the based homotopy class of $\gamma_1$ defines an element in $\pi_1N$. We consider the covering space corresponding to the subgroup generated by this element. The lifts of $\gamma_1$ and $\gamma_2$ to this cover are again simple closed curves. The singular $1$-form $d\theta$ and the vector field $X$ also lift as globally defined objects; this justifies the claim about the extension of $d\theta$ to $\gamma_2$. Up to a change of sign of $d\theta$, we can assume that at any point of $\gamma_1$, the tangent vector $v$ satisfies $d\theta(v)>0$, and similarly for $\gamma_2$.

We will now work entirely on the covering space. Note that the zeroes of $X$ are either sources or sinks, i.e.~they all have index $0$. The obstruction against finding a homotopy from $\gamma_1$ to $\gamma_2$ transverse to $\ker d\theta$ is the possible existence of zeroes of $X$, over which we cannot pull $\gamma_2$ without losing transversality. 
In a first step we can make $\gamma_1$ and $\gamma_2$ transverse to one another, keeping transversality to $\ker d\theta$. If $\gamma_1$ and $\gamma_2$ do not intersect, then they bound an annulus and are both transverse to the vector field $X$. Poincar\'e-Hopf now shows that there cannot be any zeroes in the interior of the annulus, hence we can construct a transverse homotopy.

To try to reduce to this case, we can divide both curves into subarcs in such a way that any subarc of $\gamma_1$ together with a corresponding subarc of $\gamma_2$ bounds a disk. If such a disk does not contain a zero of $X$, we can slide this subarc of $\gamma_2$ over the disk to remove two intersection points. After these modifications two possibilities remain: either $\gamma_1$ and $\gamma_2$ do not intersect anymore and we are done, or they still do. In this case, again Poincar\'e-Hopf shows that any  disk bounded by two corresponding subarcs of $\gamma_1$ and $\gamma_2$ contain exactly one zero of $X$. Now direct inspection shows that on one of the subarc the tangent vector must satisfy $d\theta(v)<0$, which we had excluded. Hence we are left with the first case and the proof is finished.  
\end{pf}

The fact that there can exist homotopic but non-equivalent splittable curves is of practical significance for the development of the deformation theory, since it implies that specifying the homotopy class of a new meridian is not enough to determine a splitting. It is actually not difficult to construct a spherical cone-surface with smooth part $N$ and a homotopy class on $N$ containing two non-equivalent splittable curves, however all the examples we have found so far have some large cone angles (i.e.~larger than $2\pi$). Thus it remains an open question whether such examples can exist in our setting, where all the cone angles are smaller than $2\pi$.
Of course, intuition suggests that many homotopy classes will not contain any splittable curve at all; but as we will see later, infinitely many non-homotopic curves satisfying the splitting condition can exist on a spherical cone-surface.

%

\medskip

A related, useful notion is that of the cone angle along a curve:
\begin{Def}
Let $N$ be an incomplete oriented spherical surface. We say that a closed curve on $N$ is non-singular if its holonomy is not trivial and the $1$-form $d\theta$ along it is not singular. The cone angle along a non-singular curve $\gamma$ is then defined as $| \int_\gamma d\theta|$.
\end{Def}

It is an easy exercise to check that if $\gamma$ is a loop around a conic point of a cone-surface, then the cone angle along $\gamma$ is exactly the cone angle of this conic point; moreover, the cone angle along a curve is always equal modulo $2\pi$ to (plus or minus) the rotation angle of its holonomy. Note also that the cone angle along $\gamma$ may very well be greater than $2\pi$, even if $N$ is the smooth part of a spherical cone-surface with cone angles smaller than $2\pi$ and $\gamma$ is a simple closed curve. Finally, we emphasize that this quantity is \emph{not} homotopy-invariant: if $H$ is a homotopy from $\gamma$ to $\gamma'$ that crosses a singular point of $d\theta$, then one easily shows that $\int_\gamma d\theta - \int_{\gamma'} d\theta = \pm 2\pi$, where the exact sign depends on the direction of the crossing. In particular, if $\gamma$ is such that $d\theta$ is singular along it, then we can always find a small perturbation of $\gamma$ for which this is no longer the case, but the resulting cone angle may depend on the chosen perturbation. Of course, since this kind of homotopy is forbidden for equivalent splittable curves, the notion of cone angle is well-defined for an equivalence class of splittable curves. As we will see in section \ref{sec:split}, this quantity is actually equal to the cone angle of the new edge obtained by splitting a vertex along a splittable curve.

\subsubsection*{Some examples}\label{sec:examples}

\begin{figure}
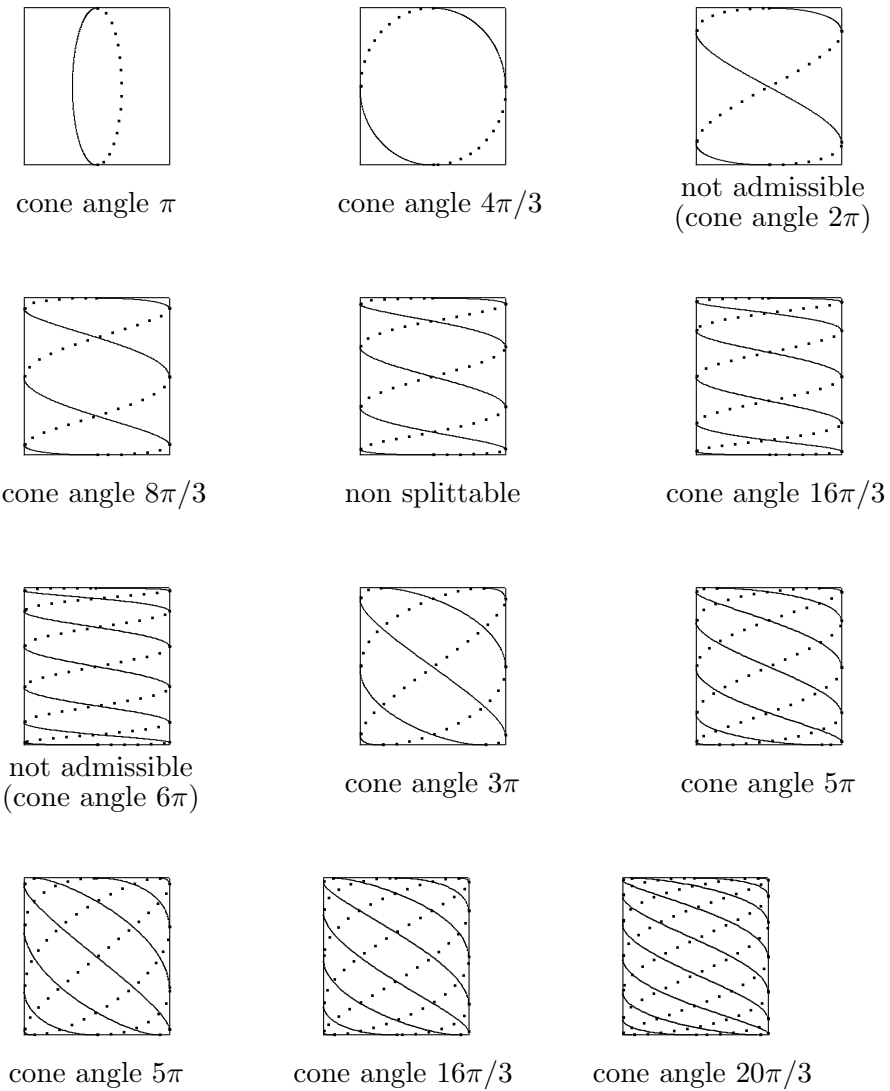

\include{Lissajous}
\caption{The cone angle along some simple closed curves on $S$.}
\label{fig:Lissajous}
\end{figure}

Let $S$ be the spherical cone-surface obtained as the double of a spherical square with angle $2\pi/3$. This cone-surface is particularly easy to study: it has a non-trivial isometry group, and more importantly, the image of its holonomy representation is finite in $\SU(2)$ (it covers the group of direct isometries of a tetrahedron). It turns out that on $S$, there exist infinitely many non-homotopic splittable curves, and that the cone angles along such curves can grow arbitrarily large. In Fig.~\ref{fig:Lissajous} we give some of the simplest curves on $S$; except for three of them, they are all splittable (actually even the two non-admissible curves satisfy an adapted splitting condition, see the end of section \ref{sec:split}).

The first seven depicted curves are obtained by applying iteratively the same half Dehn twist. Let $\gamma_k$ be the curve obtained after $k$ such half Dehn twists. Then one can verify that the homotopy class of $\gamma_k$ always contains a splittable curve (unique up to equivalence) when $k$ is odd, and that its cone angle is equal to $k\pi+\pi/3$ if $k \equiv 1 \mod 4$, and $k\pi-\pi/3$ if $k \equiv 3 \mod 4$. In particular, it goes to infinity as $k$ grows to infinity.

\subsubsection*{Bounding the cone angle}

We have seen that even if a curve $\gamma$ satisfies the splitting condition, the cone angle along $\gamma$ may well be larger than $2\pi$. Hence, in general, a splitting deformation will not preserve the class of hyperbolic cone-$3$-manifolds with cone angles less than $2\pi$. However, if $\gamma$ has special extrinsic geometry, then the cone angle along $\gamma$ can in fact
be bounded from above:

\begin{prop}\label{bounding_length}
Let $L$ be a spherical cone-surface with cone angles less than $2\pi$ and $N$
its smooth part. If $\gamma$ is a simple closed curve on $N$ of constant
geodesic curvature $\kappa \in [0,\infty)$, then the length of $\gamma$ is
less than or equal to $2\pi / \sqrt{1 + \kappa^2}$, which is precisely the length of such a curve on the standard round 2-sphere.
\end{prop}

\begin{pf}
It is a classical result of Toponogov, cf.~\cite{Top}, that a simple closed geodesic on a simply
connected smooth surface with Gauss curvature $K \geq 1$ has length less
than or equal to $2\pi$. Moreover, equality can only occur on the standard
round 2-sphere. An exposition of this appears in \cite{Kli} on p.~297. Applying this
result to smoothings of $L$ yields the same statement for simple closed geodesics on $N$. 

Having established the assertion for simple closed geodesics, we can argue
as follows in the general case: Suppose there exists a spherical cone-surface with
cone angles less than $2\pi$ and a simple closed curve $\gamma$ of constant
geodesic curvature $\kappa$, but with length $l> 2\pi/\sqrt{1 + \kappa^2}$. Clearly this curve cannot be
peripheral, since the cone angles are less than $2\pi$. If we cut our cone-surface along $\gamma$, we
obtain a piece with convex boundary and
a piece with concave boundary. Take the convex piece and add the region of a spherical
football of cone angle $\alpha=l\,\sqrt{1 + \kappa^2}$ (i.e.~the spherical suspension over a circle of length $\alpha$) bounded by the closed
geodesic of length $\alpha>2\pi$ and the circle of constant geodesic curvature $\kappa$. We may now
double to obtain a spherical
cone-surface with cone angles less than $2\pi$ containing a simple closed
geodesic of length greater than $2\pi$, which contradicts Topogonov's result.
\end{pf}

An immediate consequence is the following:
\begin{cor}
  If $\gamma$ is as above and has non-trivial holonomy, then the cone angle along $\gamma$ is less than $2\pi$.
\end{cor}

\subsection{Splitting deformations}\label{sec:split}

\subsubsection*{The splittable case}

Let $S$ be a spherical cone-surface with smooth part $N$ and holonomy representation $\varrho : \pi_1 N \to K \subset \SL_2(\C)$, and let $\nu$ be a simple closed curve on $N$ with non-trivial holonomy. Let $C$ be the untruncated hyperbolic cone over $S$. It is a complete, infinite volume hyperbolic cone-$3$-manifold, with a unique singular vertex (the summit of the cone) and singular edges corresponding to the cone points of $S$. We want to construct a one-parameter family of cone-manifolds, obtained from $C$ by splitting its vertex in such a way that the new edge has $\nu$ as a meridian. If $\nu$ satisfies the splitting condition, there is a natural way to proceed, and the resulting one-parameter family $C_l$, $l \geq 0$, is then called a \emph{splitting deformation of} $C$ \emph{along} $\nu$.

Let $\alpha = |\int_\nu d\theta|$ be the cone angle of $\nu$. The football of angle $\alpha$, denoted by $\mathbb{S}^2_\alpha$, is defined as the only (up to isometry) spherical cone-surface with underlying space $\mathbb{S}^2$ and two cone-points of angle $\alpha$; it can be constructed as the spherical suspension over a circle of length $\alpha$. More relevantly, in our case, if $r$ and $\theta$ are polar coordinates on the universal cover of $\mathbb{S}^2\setminus Fix(\varrho(\nu))$, then $\mathbb{S}^2_\alpha$ is the metric completion of the quotient of this universal cover under the identification $(r,\theta) \sim (r,\theta+\alpha)$. It is then possible to embed isometrically a neighbourhood of $\nu$ in $\mathbb{S}^2_\alpha$; the image of $\nu$ by this embedding will be denoted by $\nu'$.

The hyperbolic cone $C(\nu)$ over $\nu$ separates $C$ into two ``cone-manifolds with boundary'' $C_1$ and $C_2$; if $\nu$ is smooth, the boundary $\partial C_i = C(\nu)$ is a smooth surface except at the singular vertex. Now let $T$ be the untruncated hyperbolic cone over $\mathbb{S}^2_\alpha$. It is an infinite singular tube, with a unique singular edge of cone angle $\alpha$. The cone over $\nu'$ is a surface $B_1 = C(\nu')$ in $T$, and since $\nu$ and $\nu'$ have isometric neighbourhoods, $C(\nu)$ and $C(\nu')$ are isometric and have isometric conic neighbourhoods. Let $B_2$ be the image of $B_1 = C(\nu')$ by a hyperbolic translation of length 
$l$ along the singular edge of $T$. The fact that $\nu$ satisfies the splitting condition implies that $B_1$ and $B_2$ do not intersect; in particular, they bound a region $U$. We can now define the new (complete, infinite volume) cone-manifold $C_l$ as $$C_l =  C_1 \cup U \cup C_2 / \sim,$$
 where the boundaries $\partial C_i$ and $B_i$ are identified; by construction, the meridian of the new edge is homotopic to $\nu$ and its cone angle is equal to the one along $\nu$. 

This method also allows us to split $C$ along a family of disjoint splittable curves. More precisely, if $\nu_1,\ldots, \nu_{p}$ are pairwise disjoint, non-homotopic curves on $N$ satisfying the splitting condition, then we can split $C$ along this family by proceeding as above: we first cut $C$ along the $C(\nu_i)$ (which only intersect at the singular vertex), then glue parts of singular tubes between the pieces. By restricting this construction, we can apply it to a neighbourhood of a singular vertex $v$ of $X$; in fact, we can thus produce a splitting deformation of the neighbourhood $U_\varepsilon(\Sigma)$ of the singular locus. 

\begin{figure}
\centering 
\includegraphics{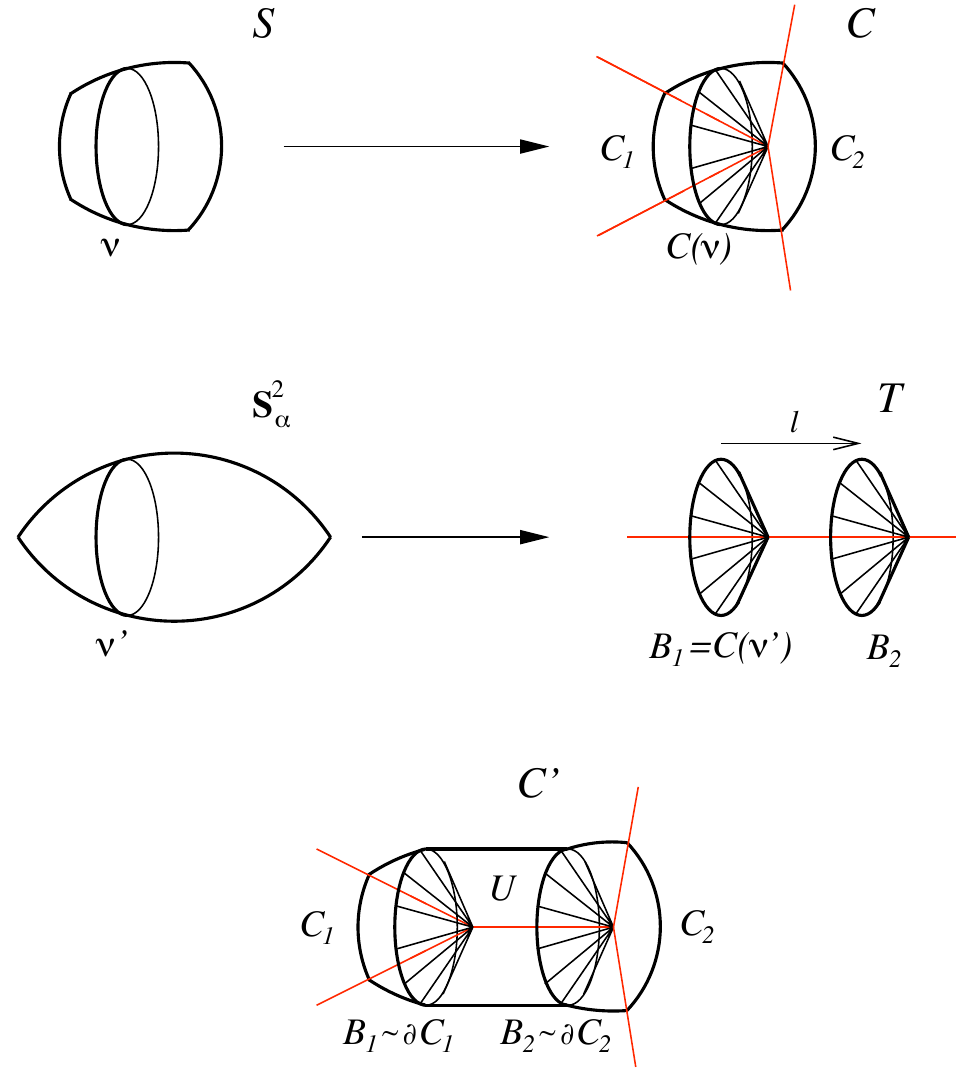}
\caption{Construction of model splittings}
\end{figure}

\begin{rem} If we replace $\nu$ by an equivalent splittable curve in this construction, then it is clear that the resulting cone-manifolds are the same. However, replacing $\nu$ by a homotopic but non-equivalent splittable curve will yield a different splitting deformation $C'_l$ of $C$, even though the meridian of the new edge is still homotopic to $\nu$. The two families $C_l$ and $C'_l$ should be thought as splitting $C$ into opposite directions, and we will see that they can be distinguished by their holonomy characters.
\end{rem}

\subsubsection*{More general splitting deformations}

The splitting condition that we introduce in section \ref{sec:splitcond} is rather simple to formulate and to verify; however, it is not strictly necessary for the existence of splitting deformations. In fact, the above construction can be carried out for a slightly larger class of curves, either non-splittable or with trivial holonomy.

We first consider the case of a non-peripheral, non-singular simple closed curve $\nu$ on the smooth part $N$ of a spherical cone-surface $S$. The curve $\nu$ separates $S$ into two parts $S_1$ and $S_2$, and we want to construct the splitting deformation of the hyperbolic cone $C(S)$ by gluing a part of a singular tube between $C_1=C(S_1)$ and $C_2=C(S_2)$. We know that a neighbourhood of $\nu$ can be isometrically immersed in the football $\mathbb{S}^2_\alpha$; we will furthermore assume that this immersion is an embedding. 

\begin{figure}
\centering
\includegraphics[scale=0.65]{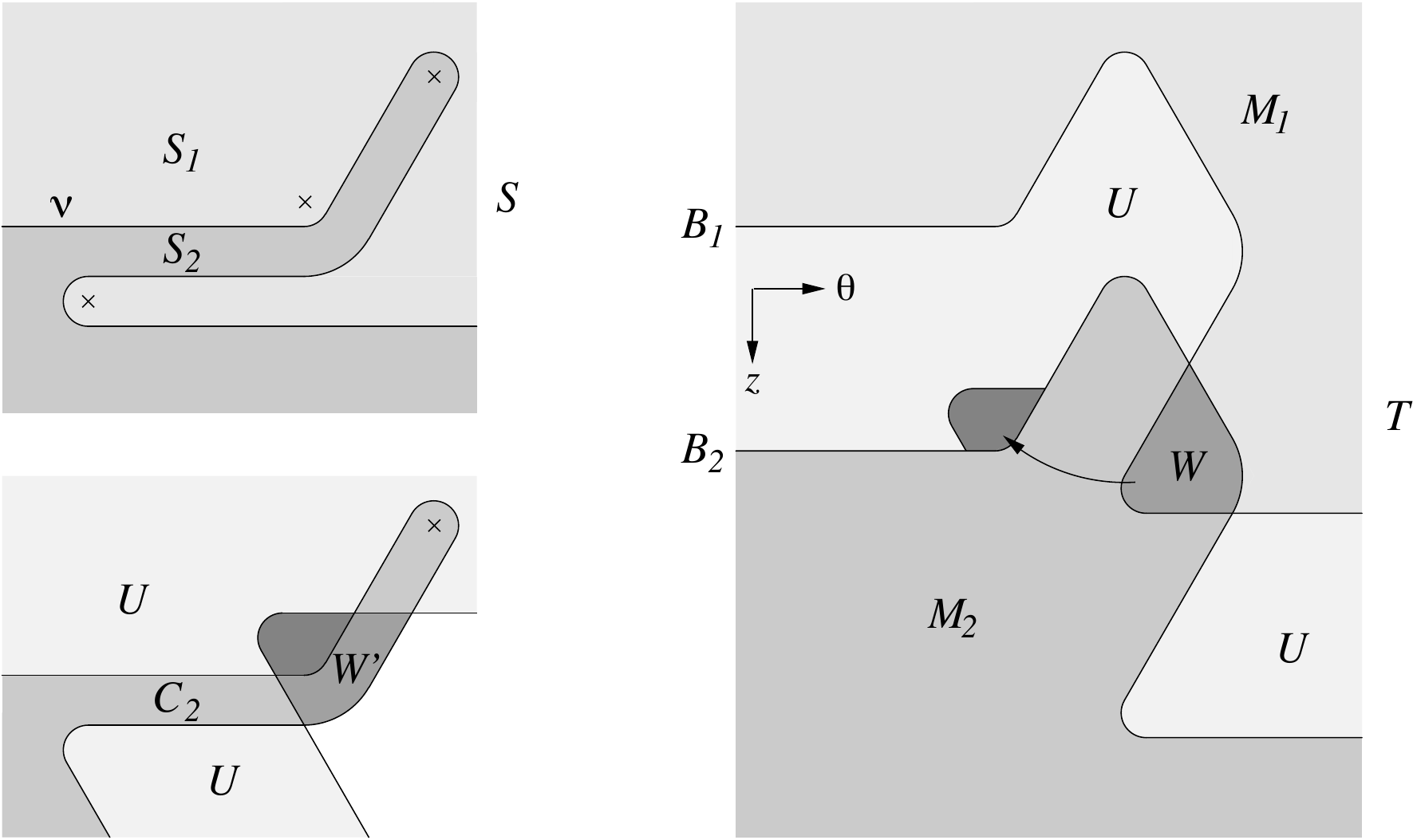}
\caption{A generalized splitting deformation}
\label{fig:intersplit}
\end{figure}

As before, let $\nu'$ be the image of $\nu$ by this embedding, and let $B_1$ be the cone over $\nu'$ in $T=C(\mathbb{S}^2_\alpha)$. The surface $B_2$ is obtained as the image of $B_1$ by a hyperbolic translation of length $l$ along the singular axis of $T$. Now $B_1$ (resp.~$B_2$) separates $T$ into two half-spaces; let $M_1$ (resp.~$M_2$) be the one lying in the negative (resp.~positive) direction with respect to the direction of the hyperbolic translation above. Several cases can happen:

\begin{itemize}
\item The two half-spaces $M_1$ and $M_2$ (and consequently also the two surfaces $B_1$ and $B_2$) do not intersect. This is the ``good'' case, happening in particular when $\nu$ is splittable. Then we can proceed as before: $B_1$ and $B_2$ bound a region $U= T \setminus (M_1 \cup M_2)$, and the deformed cone-manifold is constructed as $C_1 \cup U \cup C_2$ with boundaries identified.

\item The two half-spaces $M_1$ and $M_2$ intersect. Let $U=T \setminus (M_1 \cup M_2)$ and $W= M_1\cap M_2$; for the simplicity of the discussion we will assume that $W$ is connected, but what follows can be applied to all the components of $W$. The boundary of $W$ is composed of two parts, $\partial W_1 \subset B_1$ and $\partial W_2 \subset B_2$. The isometry between neighbourhoods of $\nu$ and $\nu'$ yields an isometric embedding of a (conic) neighbourhood of $B_1$ in $M_1$ into $C_1$; starting from a neighbourhood of $\partial W_1$, we can try to extend this to an isometric embedding of $W$ into $C_1$. We can proceed similarly for embedding isometrically $W$ into $C_2$. Here again, several cases can happen:

\begin{enumerate}
\item The region $W$ can be embedded isometrically as $W'$ in e.g.~$C_2$. Then the deformed cone-manifold can be constructed as $C_1 \cup U \cup (C_2 \setminus W')$, with the consistent boundary identifications. This case arises for instance when $\nu$ could have been homotoped into a splittable curve.

\item The region $W$ cannot be embedded in $C_1$ or $C_2$, but can be embedded isometrically as $W'$ in $C_2\cup U$, where the boundaries $\partial C_2$ and $B_2$ are identified (the case is of course similar if $W$ embeds isometrically in $C_1\cup U$). Then the deformed cone-manifold can be constructed as $C_1 \cup ((C_2 \cup U) \setminus W')$, with the consistent boundary identifications. We give in Fig.~\ref{fig:intersplit} an example of this construction: in the top-left corner is represented a part of a spherical cone-surface $S$ with three of its cone points, and the (non-splittable) simple closed curve $\nu$ that separates it into $S_1$ and $S_2$. The right side of the picture is a sketch of $T$ and its various regions. If we use cylindrical coordinates $(r,z,\theta)$ in $T$ (as in the introduction), then this sketch should be thought as a slice $\{r=\mbox{constant}\}$; the foliation given by $\ker d\theta$ is vertical, and we see at once that $\nu$ is not splittable. The bottom-left corner depicts (a slice of) $C_2\cup U$, as well as $W'$, the isometric image of $W$. The darkest part of $W'$ lies in $U$ and is also represented on the right.

\item the region $W$ cannot be embedded isometrically in $C_1 \cup U$ or $C_2 \cup U$, for instance because such an embedding would encounter a singular edge. Then the construction of the deformed cone-manifold along these lines fails.

\end{enumerate}
\end{itemize}

\begin{figure}
\centering
\includegraphics[scale=0.65]{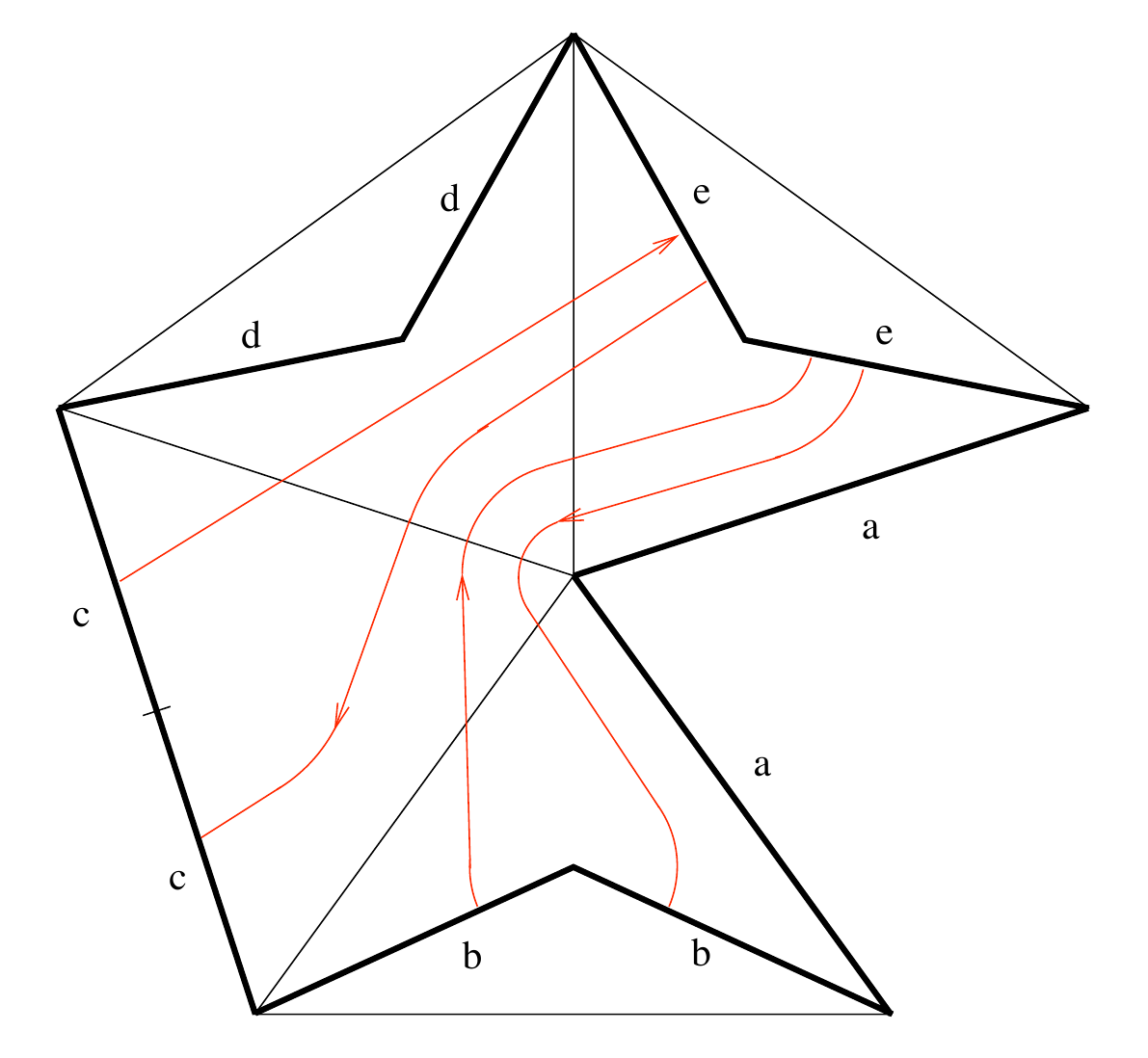}
\caption{A generalized splittable curve}
\label{fig:counterex}
\end{figure}

The case 2 is of course the most interesting: the curve $\nu$ is then called a \emph{generalized splittable curve}. Fig.~\ref{fig:counterex} shows an example where this situation arises. It represents a spherical cone-surface; the thick boundary edges should be identified according to their labels. All the big triangles are equilateral and have angles $2\pi/5$, so that the holonomy representation actually has values in the icosahedron symmetry group. It has $5$ cone points of cone angles $\pi$, $4\pi/3$ (three times) and $8\pi/5$. The cone angle along the simple closed curve $\nu$ depicted in the picture is $8\pi/5$. The homotopy class of $\nu$ is not splittable; however, splitting deformations along this curve are possible, according to the above construction. 

Splitting deformations can also be defined for curves with trivial holonomy (although they never can be part of an admissible pair-of-pants decomposition of $\partial \bar M_\varepsilon$). 
Let $\nu$ be a non-peripheral, simple closed curve with trivial holonomy on $N$. Let $\tilde \nu$ be a lift of $\nu$ to the universal cover $\tilde N$. Since the curve $\nu$ has trivial holonomy, the image of $\tilde \nu$ under the developing map is a closed curve on $\mathbb{S}^2$. Now every point $p$ of $\mathbb{S}^2$ defines polar coordinates $(r,\theta)$ on $\mathbb{S}^2 \setminus \{p,-p\}$, and the $1$-form $d\theta$ descends to a well-defined $1$-form, still denoted $d\theta$, in a neighbourhood of $\nu$ (non-singular if the image of $\tilde \nu$ does not intersect $\{p,-p\}$). We can then proceed exactly as previously. In particular, we can define a splitting condition for $\nu$ by requiring that it is transverse to such a $1$-form $d\theta$ corresponding to some point $p\in \mathbb{S}^2$; similarly, the cone angle along $\nu$ with respect to $d\theta$ is simply defined as $|\int_{\nu} d\theta|$, and is necessarily a multiple of $2\pi$.

A splitting deformation obtained in that way (provided that it can be constructed, either because $\nu$ is splittable or that the generalized construction works) depends obviously of the choice of $p$, and it is quite unclear whether it can be extended to the whole of $M$. Note that if a new edge $e$ has cone angle $2\pi$, then the Riemannian metric is actually not singular along $e$, and this ``removable edge'' can be included in the smooth part. This operation, however, changes the topological type of the smooth part, and thus does not fit well within our $\Defo(M)$-framework.

\subsection{Deformations and the holonomy character}\label{sec:defandhol}

The splitting deformations of the singular tube $U_\varepsilon(\Sigma)$ that
we have just constructed are closely related to the complex Hamiltonian (twist) flows introduced in section \ref{sec:complex-sympl}. Since $U_\varepsilon(\Sigma)$ retracts by deformation on $\partial \bar M_\varepsilon$, they have the same character variety, and a hyperbolic metric on the singular tube yields a holonomy character $\chi \in X(\pi_1 U_\varepsilon(\Sigma), \SL_2(\C)) = X(\pi_1 \partial \bar M_\varepsilon, \SL_2(\C))$. Let $\chi_0=[\rho_0]$ be the character induced by the hyperbolic structure of the cone-manifold $X$. 

\begin{lemma}\label{lem:sign}
Let $N$ be the smooth part of a link's vertex and let $\nu$ be a splittable
curve on $N$. There exists a sign $s\in \{-1,1\}$ such that the deformation of
$U_\varepsilon(\Sigma)$ obtained by splitting along $\nu$ has holonomy
character $\chi' = \phi_{ist}(\chi_0)$, where $\phi$ is the Hamiltonian flow associated to $\tr_{\nu}$,  $t=l\left(4-\tr_{\nu}(\chi_0)^2\right)^{-1/2}$, and $l$ is the length of the new edge. 
\end{lemma}

\begin{pf}
We begin with the construction of the sign $s$. Let $d\theta$ be the standard $1$-form along $\nu$. Up to a change of sign of $\theta$, we can assume that the integral of $d\theta$ along $\nu$ is positive (hence equal to the cone angle $\alpha$ along $\nu$). The vector field $\partial / \partial \theta$ is a Killing vector field in a neighbourhood of $\nu$, and as such it corresponds to a parallel section $\sigma_{\partial / \partial \theta}$ of the $\sl_2(\C)$-bundle $\E_{\rho_0}$. 
Up to conjugacy, we can assume that the restriction of $\rho_0$ to $N$ has values in $\SU(2)$, 
and that the fibers of $\E_{\rho_0}$ over $\nu$ are identified with $\sl_2(\C)$ in such a way that $\sigma_{\partial/\partial \theta}$ is the constant section $\begin{bmatrix} i/2 & 0 \\ 0 & -i/2\end{bmatrix}$. This implies the more precise statement that 
$$\rho_0(\nu) = c \begin{bmatrix} \exp(i\alpha /2) & 0 \\ 0 & \exp(-i\alpha /2) \end{bmatrix},$$ 
for a sign $c \in \{-1,1\}$ (this sign depends of the choice of a lifting from $\Isom(\H^3) \simeq \PSL_2(\C)$ to $\SL_2(\C)$). Since the holonomy of $\nu$ is not trivial, $\alpha \neq 0 \mbox{ mod } 2\pi$, and we can let 
$$s = - c \mathop{\rm sign}(\sin(\alpha/2));$$ 
it is equal to minus the sign of the imaginary part of the first diagonal coefficient of $\rho_0(\nu)$, once the section $\sigma_{\partial/\partial \theta}$ is correctly identified. (Note that this identification is only necessary when $\alpha = \pi \mbox{ mod } 2\pi$; otherwise it is sufficient to require that $\rho_0(\nu)$ has the above expression.)

The curve $\nu$ separates $N$ in two subsurfaces; let  $N_1$ be the one lying to the left of $\nu$, and $N_2$ the other one. We can choose a representative $\rho'$ of $\chi'$ such that the restrictions of $\rho_0$ and $\rho'$ to $N_1$ coincide (and in particular have values in $\SU(2)$). 
Now let $\zeta = \begin{bmatrix} \exp(l/2) & 0 \\ 0 & \exp(-l/2) \end{bmatrix}$; this is a hyperbolic isometry of translation length $l$. Then it is easy to check that on $\pi_1(N_1 \cup N_2)$, $\rho'$ and $\rho_0$ are related by the formula
$$
\rho'(\gamma) = \left\{ 
\begin{array}{c@{\quad:\quad}c}
\rho_0(\gamma) & \gamma \in \pi_1N_1\\ 
\zeta \rho_0(\gamma)\zeta^{-1} & \gamma \in \pi_1N_2
\end{array}
\right.
$$
This expression matches the description of the complex twist flow from section \ref{sec:complex-sympl}. More precisely, if $\psi$ denotes the complex twist flow associated to the function $\tr_{\nu}$, then $\rho'$ is equal to $\psi_z(\rho_0)$ for any $z$ such that $\zeta = \pm\exp(zF(\rho_0(\nu)))$, which is equivalent to $\exp(l/2) = \pm\exp(izc\sin(\alpha/2))$. All such values of $z$ have the same imaginary part, and a solution is given by 
$$z= -icl(2\sin(\alpha/2))^{-1} = isl(2\sqrt{1-\cos(\alpha/2)})^{-1} = isl(4-\tr_\nu(\chi_0))^{-1/2}.$$
In particular, if we let $t=l(4-\tr_\nu(\chi_0))^{-1/2}$, we obtain that the
diffeomorphism $\phi_{ist}$ associated to the Hamiltonian flow of $\tr_{\nu}$ maps $\chi_0$ to $\chi'$.
\end{pf}

\begin{rem} We have seen that there can exist on $N$ a splittable curve $\nu'$, homotopic but not equivalent to $\nu$. Then one can check that the signs associated to $\nu$ and $\nu'$ are opposite, so that in terms of holonomy, splitting along $\nu$ or along $\nu'$ corresponds to following the Hamiltonian flow $\phi_{it}$ in two different directions. 
\end{rem}

We can now deal with the remaining deformations of $U_\varepsilon(\Sigma)$, and relate them to the action-angle coordinates introduced in section \ref{sec:actionangle}. The following proposition is the main result of this section:

\begin{prop}\label{prop:defoU} 
Let $\vec{\mu}$ be the meridian set of $U_\varepsilon(\Sigma)$, and let $\vec{\nu}$ be a family of curves on $\coprod_j N_j$ such that the family $\mathcal{C}=\vec{\mu}\cup\vec{\nu}$ is up to homotopy an
admissible pair-of-pants decomposition of $\partial \bar M_\varepsilon$. In particular, the map $(\tr_{\vec{\mu}},\tr_{\vec{\nu}},\tau_{\vec{\mu}},\tau_{\vec{\nu}})$ is a system of local coordinates on $X_{irr}(\pi_1 \partial \bar M_\varepsilon, \SL_2(\C))$. Let $\chi'$ be a character close to $\chi_0$, such that $(\tr_{\vec \mu}(\chi'),\tr_{\vec \nu}(\chi')) \in (-2,2)^{3g-3}$. Assume furthermore that the curves in $\vec{\nu}$ for which $\Im \tau_{\nu_i}(\chi') \neq 0$ are splittable, disjoint, and with associated signs $s_i$ such that $s_i \Im \tau_{\nu_i}(\chi') > 0$. Then there exists a deformation of $U_\varepsilon(\Sigma)$, whose holonomy character is precisely $\chi'$.
\end{prop}

\begin{pf}
We begin by constructing a deformation $U'$ of $U_\varepsilon(\Sigma)$ whose holonomy character $\chi'_0$ satisfies $(\tr_{\vec{\mu}},\tr_{\vec{\nu}}) (\chi'_0)= (\tr_{\vec{\mu}},\tr_{\vec{\nu}}) (\chi')$, $\tau_{\vec{\mu}}(\chi'_0)= 0$, and $\tau_{\vec{\nu}} (\chi'_0) = \Re(\tau_{\vec{\nu}}(\chi'))$. Due to our choice of action-angle coordinates, this means that the restriction to any vertex's link $N_j$ of a representative $\rho'_0$ of $\chi'_0$ has values in a maximal compact subgroup of $\SL_2(\C)$; this is of course also true for $\rho_0$. In particular, we can consider the restrictions of $\chi_0$ and $\chi'_0$ as elements of $X(\pi_1 N_j, \SU(2))$. Now the local homeomorphism (akin to \eqref{eq:hol}) between the space $\Defo(N_j)$ of spherical structures on $N_j$ and $X(\pi_1 N_j, \SU(2))$ means that we can deform the link $L_j$ into a spherical cone-surface $L'_j$ whose holonomy character is given by $\chi'_0$. Taking the hyperbolic cone over $L'_j$, we obtain the corresponding deformation of a neighbourhood of $v_j$, and then the deformation $U'$ of $U_\varepsilon(\Sigma)$ associated to $\chi'_0$. A second, easy step is to construct a deformation $U''$ whose holonomy character $\chi''_0$ is such that 
$(\tr_{\vec{\mu}},\tr_{\vec{\nu}},\tau_{\vec{\mu}},\tau_{\vec{\nu}}) (\chi''_0)= (\tr_{\vec{\mu}}(\chi'),\tr_{\vec{\nu}}(\chi'),\tau_{\vec{\mu}}(\chi'),\Re(\tau_{\vec{\nu}}(\chi'))$. For this, all we have to do is to change the lengths and/or twist parameters of the existing edges, see \cite{Maz-GM} section 3 for more details. 

If $\chi'$ and$\chi_0$ are sufficiently close, Lemma \ref{lem:open} implies that the curves $\nu_i$ for which $\Im \tau_{\nu_i}(\chi') \neq 0$ are still splittable and disjoint in $U''$; furthermore, it is easy to check that the associated signs $s_i$ given by Lemma \ref{lem:sign} are preserved if the deformation is small enough. We can thus apply the splitting deformations described above simultaneously to all the curves $\nu_i$ for which $\Im \tau_{\nu_i}(\chi') \neq 0$, inserting new edges of lengths $l_i=s_i \Im(\tau_{\nu_i}(\chi')) \sqrt{4-\tr_{\nu_i}(\chi')^2}$. Lemma \ref{lem:sign} then shows that the resulting deformation of $U_\varepsilon(\Sigma)$ has the correct holonomy character $\chi'$.
\end{pf}

The above proposition shows the existence of the deformation for the character $\chi'$, but not the uniqueness. Actually, in the correspondence \eqref{eq:hol}, the deformation of the hyperbolic structure on $U_\varepsilon(\Sigma)$ corresponding to a character close to $\chi_0$ is unique only up to the relation induced by thickening (and isotopy), see \cite{CHK}. So $\chi'$ determines the hyperbolic metric only in a compact core, and it can happen that this metric admits different completions as a deformation of the cone-manifold structure on $U_\varepsilon(\Sigma)$. The following definition lifts this ambiguity:

\begin{Def}\label{def:nucomp} 
A deformation of the hyperbolic cone-manifold structure on $U_\varepsilon(\Sigma)$ is called $\vec \nu$-compatible if it is obtained as a sequence of deformations as in the proof of Proposition \ref{prop:defoU}.
\end{Def}

Note that by construction, a $\vec{\nu}$-compatible deformation is not split along any non-splittable curve; furthermore if $\vec{\nu}$ contains two non-disjoint curves then the deformation is split along only one of them.

\section{The local shape of the deformation space}\label{sec:localshape}

Having completed our study of the deformations of the singular tube, we now turn to the deformations of cone-manifolds. We actually limit ourselves to compatible deformations, as defined below. Using the results of sections \ref{sec:coord} and \ref{sec:defosingtube}, this enables us to determine the possible deformations of the holonomy representation and of the geometric structure in a neighbourhood of the singular locus. With these two ingredients we can then proceed to prove our main results (Theorems \ref{thm:construction} and \ref{thm:main}). As an application of these results, we explain the structure of the deformation space of $X$ and detail what happens for doubles of polyhedra.

\subsection{Compatible cone-manifold structures}

We classify the deformations of $X$ according to the deformations they induce on the singular tube $U_\varepsilon(\Sigma)$.

\begin{Def}\label{def:compatible} Let $\vec{\mu}$ be the meridian set of $X$, and let $\vec{\nu}$ be a family of curves on $\coprod_j N_j$ such that the family $\mathcal{C}=\vec{\mu}\cup\vec{\nu}$ is up to homotopy an admissible pair-of-pants decomposition of $\partial \bar M_\varepsilon$.
A deformation of the hyperbolic cone-manifold structure on $X$ is called $\vec \nu$-compatible if its restriction to $U_\varepsilon(\Sigma)$ is $\vec \nu$-compatible. We denote by $C_{-1}(X,\vec \nu)$ be the space of $\vec \nu$-compatible deformations of $X$.
\end{Def}

We identify an element in $C_{-1}(X,\vec \nu)$ with the corresponding deformation of the smooth part, i.e.~with an element in $\Defo(M)$. Hence we may topologize $C_{-1}(X,\vec \nu)$ as a subset of $\Defo(M)$. 

If $X$ has cone angles less than $2\pi$, then $\mathcal{C}$ defines a local coordinate system on $X_{irr}(\pi_1M,\SL_2(\C))$ near the holonomy character $\chi_0 = [\hol]$ of $X$. This coordinate system is given by Theorem \ref{thm:param}, which asserts that the map 
\begin{align*}
\Phi_{\mathcal{C}} : \mathcal{U} \cap X_{irr}(\pi_1 M,\SL_2(\C)) & \rightarrow \C^N \times \R^{3g-3-N} \times \R^{3g-3-N} \\
\chi & \mapsto ( \tr_{\vec \mu}(\chi), \Im \tau_{\vec \nu}(\chi), \Im \tr_{\vec \nu}(\chi))
\end{align*}
is a local diffeomorphism at $\chi_0$. Since via \eqref{eq:hol} the deformation space $\Defo(M)$ is locally homeomorphic to the character variety $X_{irr}(\pi_1M,\SL_2(\C))$, we also obtain in this way a local coordinate chart on $\Defo(M)$.
There is of course a connection between this coordinate map $\Phi_{\mathcal{C}}$ and the subspace of $\vec \nu$-compatible deformations of $X$. By definition of the action-angle coordinates, we know that $\Phi_{\mathcal{C}}(\chi_0) \in ({-}2,2)^N \times \{0\}^{6g-6-2N}$. Now if a character $\chi' \in \mathcal{U}$ corresponds to an element $X' \in C_{-1}(X,\vec \nu)$, then every curve in $\mathcal{C}$ is either homotopic to the the meridian of an edge or is up to homotopy contained in the link of a vertex of $X'$; in any case, its holonomy is elliptic. This means that $\Phi_{\mathcal{C}}$ restricted to $\mathcal{U} \cap C_{-1}(X,\vec \nu)$ actually has values in $({-}2,2)^N \times \R^{3g-3-N} \times \{0\}^{3g-3-N}$. Using Lemma \ref{lem:sign}, we can actually be more precise, since we know that $\Im \tau_{\vec{\nu}_i}(\chi')$ must have sign $s_i$ if $\nu_i$ is splittable, and vanish otherwise. 

In the same lemma, we have also seen a formula for the length of a new edge:
$$\ell_i =  s_i \Im \tau_{\nu_i}(\chi') \sqrt{4-\tr_{\nu_i}(\chi')^2}.$$
This suggests that on $C_{-1}(X,\vec\nu)$, we can use more geometrically meaningful coordinates. For simplicity, we assume that the family $\vec{\nu}$ is sorted so that the curves $\nu_1$ to $\nu_{n_1}$ are splittable and the others are not.  Let $\R_s^{n_1}$ be the subset of tuples $(x_1,\ldots,x_{n_1})$ such that $s_ix_i\geq 0$.
Now for each curve $\mu_i \in \vec{\mu}$, there exists a sign $c_i \in \{-1,1\}$ such that $\tr_{\mu_i}(\chi_0) = 2 c_i \cos(\alpha_i/2)$ (this is the same definition as in the proof of Lemma \ref{lem:sign}). Let $f_i(x) = 2 c_i \cos(x/2)$; it is a diffeomorphism from $(0,2\pi)$ to $(-2,2)$, and $f_i(\alpha_i) = \tr_{\mu_i}(\chi_0)$.

\begin{prop}\label{prop:param3}
Let $E_{\vec{\nu}} \subset \mathcal U$ denote the submanifold-with-corners  $$\Phi_{\mathcal{C}}^{-1}( \mathcal V \cap (-2,2)^N\times \R_s^{n1}\times\{0\}^{3g-3-N-n_1}\times\{0\}^{3g-3-N}).$$ 
Then the map $\Psi_{\vec{\nu}} : E_{\vec{\nu}} \to (0,2\pi)^N \times \R_{\geq 0}^{n_1}$,
\begin{multline*}
\Psi_{\vec{\nu}}(\chi) = \Big(f_1^{-1}(\tr_{\mu_1}(\chi)),\ldots,f_N^{-1}(\tr_{\mu_N}(\chi)),
 s_1 \Im \tau_{\nu_1}(\chi) \sqrt{4-\tr_{\nu_1}(\chi)^2},\\ \ldots,s_{n_1} \Im \tau_{\nu_{n_1}}(\chi) \sqrt{4-\tr_{\nu_{n_1}}(\chi)^2}\Big),
\end{multline*}
is a local diffeomorphism. Moreover, if $\chi'$ is the holonomy character of a cone-manifold structure $X' \in C_{-1}(X,\vec\nu)$ close to $X$, then 
$$\Psi_{\vec{\nu}} (\chi') = (\alpha'_1,\ldots,\alpha'_N,\ell_1,\ldots,\ell_{n_1}),$$ 
where $(\alpha'_i)_{1\leq i \leq N}$ are the cone angles of the original edges and $(\ell_i)_{1\leq i \leq n_1}$ are the lengths of the new edges.
\end{prop}

\begin{pf} It follows directly from the preceding discussion and the fact that the map $\Phi_{\mathcal{C}}$ is a local diffeomorphism (Theorem \ref{thm:param}).
\end{pf}

\subsection{Deformations of hyperbolic cone-$3$-manifolds}

We have just seen that the local diffeomorphism $\Psi_{\vec{\nu}} : E_{\vec{\nu}} \to (0,2\pi)^N \times \R_{\geq 0}^{n_1}$ induces on $C_{-1}(X,\vec\nu)$ a locally injective map, that sends an element to the vector composed of the cone angles of the original edges and the lengths of the new ones. As a final step, we now show that assuming the splitting condition, this map is locally onto, or equivalently that $C_{-1}(X,\vec\nu)$ can be locally identified with $E_{\vec{\nu}}$. 

\begin{thm}\label{thm:construction}
Let $X$ be a hyperbolic cone-$3$-manifold with meridian set $\vec{\mu}$ and cone angles $\alpha = (\alpha_1,\ldots,\alpha_N) \in (0,2\pi)^N$. Let $\vec\nu$ be a family of curves on $\coprod_{j=1}^kN_j$ such that ${\mathcal C}= \vec \mu \cup \vec \nu$ is up to homotopy an admissible pair-of-pants decomposition of $\partial \bar M_\varepsilon$. Then the $\vec{\nu}$-compatible deformation of the hyperbolic cone-manifold structure on $X$ corresponding to a vector $(\alpha',\ell) \in (0,2\pi)^N \times \R_{\geq 0}^{3g-3-N}$ close to $(\alpha,0)$ exists, if the curves $\nu_i$ with $\ell_i>0$ are splittable and disjoint.
\end{thm}

\begin{pf} 
Using Proposition \ref{prop:param3}, we know the candidate $\chi' = \Psi_{\vec{\nu}}^{-1} (\alpha',\ell)$ for the conjugacy class of the holonomy representation of $X'$. Since $\Defo(M)$ is locally homeomorphic to $X(\pi_1M,\SL_2(\C))$, we know that $\chi'$ is the holonomy character of a neighbouring hyperbolic structure on $M$. What we need to show is that this hyperbolic structure (or a thickening thereof) can be completed into a cone-manifold $X' \in C_{-1}(X,\vec{\nu})$, i.e.~such that its singular locus $\Sigma'$ is obtained from $\Sigma$ by splitting some vertices according to the curves of $\vec{\nu}$. 

By abuse of notation, the restriction of $\chi'$ to $X(\pi_1 \partial \bar M_\varepsilon, \SL_2(\C))$ will still be denoted $\chi'$. Similarly, $\chi_0$ stands for both the conjugacy class of the holonomy representation $\hol$ of $X$ and its restriction. 
Using Proposition \ref{prop:defoU}, and the fact that the curve $\nu_i$ satisfies the splitting condition and are disjoint whenever $\ell_i$ is positive, we can continuously deform the hyperbolic structure on $U_\varepsilon(\Sigma)$ into a $\vec{\nu}$-compatible one whose holonomy representation is given by $\chi'$.

Let us now recall briefly some elements of the theory of deformations of hyperbolic structure, see \cite{Goldman4} for more details. We know that a hyperbolic structure on $M$ is determined by its developing map $dev : \tilde{M} \to \H^3$ and its holonomy representation $\rho : \pi_1M \to \SL_2(\C)$. On $\tilde{M}$, we can consider the trivial bundle $\tilde{E} = \tilde{M} \times \H^3$, which admits a horizontal foliation $\mathcal{F}$ given by the constant sections. The developing map then gives rise to a section $\tilde{\sigma}$ of $\tilde{E}$, mapping $p\in \tilde{M}$ to $(p,dev(p))$; this developing section is transverse to $\mathcal{F}$. We can quotient $\tilde{E} = \tilde{M} \times \H^3$ by the action of the fundamental group of $M$ given by the holonomy representation: $(p,x) \sim (\gamma.p, \rho(\gamma)(x))$. The result is a bundle $E_\rho$ on $M$, with fiber $\H^3$, still endowed with the horizontal foliation inherited from $\mathcal{F}$. The equivariance of the developing map means that the section $\tilde{\sigma}$ descends to the developing section $\sigma$ of $E_\rho$, transverse to $\mathcal{F}$. 
What we have described here is a construction that to a hyperbolic structure on $M$ associates a triplet $(E,\mathcal{F}, \sigma)$ where $E$ is a $(\H^3,\SL_2(\C))$-bundle on $M$ with a horizontal foliation $\mathcal{F}$ and a section $\sigma$ transverse to this foliation. One can then show that any such triplet 
actually determines a hyperbolic structure on $M$.

Now let $\rho'$ be a representative of $\chi'$ close to $\hol$. Then the corresponding bundles $E_{\rho'}$ and $E_{\hol}$ are actually isomorphic, so we can identify them: $E_{\rho'} \simeq E_{\hol} \simeq E$. The two horizontal foliations $\mathcal{F}$ and $\mathcal{F'}$ are however different but close. We have seen that using the results of the previous section, we can construct on a neighbourhood $U_\varepsilon(\Sigma)$ of the singular locus of $X$ a continuous family of hyperbolic structures, joining the initial structure to one whose holonomy representation is induced by $\rho'$; this means that over $U_\varepsilon(\Sigma)$, the foliation $\mathcal{F'}$ is still transverse to $\sigma$. By compactness, this is also true over the remainder of $M$ if $\rho'$ is close enough to $\hol$, i.e.~if $(\alpha'-\alpha,\ell)$ is small enough. This implies that the triplet $(E,\mathcal{F'}, \sigma)$ determines a hyperbolic structure on $M$, whose restriction to $U_\varepsilon(\Sigma)$ is the deformation constructed above. In particular, the metric completion of $M$ is the desired cone-manifold.
\end{pf}

\begin{rem}\label{rem:2} It should be noted that the splitting assumption on the curves $\nu_i$ for which $\ell_i \neq 0$ is not indispensable; actually, the only necessary condition is that the correct deformation on $U_\varepsilon(\Sigma)$ can be constructed. In particular, we can include generalized splittable curves, as described in the last part of section \ref{sec:split}. 
\end{rem}

From this result, our main theorem follows easily. It also supposes the splitting condition, since it is easier to state, but the reader should be aware that it applies to a slightly more general situation.

\medskip

\noindent {\bf Theorem \ref{thm:main}.}
{\em Let $X$ be a hyperbolic cone-$3$-manifold with cone angles less than $2\pi$ and
meridian set $\vec \mu$. Let $\vec\nu$ be a pair-of-pants decomposition of $\coprod_{j=1}^kN_j$ such that ${\mathcal C}= \vec \mu \cup \vec \nu$ gives an admissible pair-of-pants decomposition of $\partial \bar M_\varepsilon$. If all the curves in $\vec\nu$ are splittable, then the map
\begin{equation*}
(\alpha,\ell): C_{-1}(X, \vec\nu) \rightarrow (0,2\pi)^N \times \R_{\geq 0}^{3g-3-N}
\end{equation*}
sending a $\vec\nu$-compatible cone-manifold structure to the vector composed of its original edges' cone angles and new edges' lengths, is a local homeomorphism at the given structure.}

\medskip

\begin{pf}
This follows directly from Proposition \ref{prop:param3} and Theorem \ref{thm:construction}.
\end{pf}

\subsection{Stratified structure of $C_{-1}(X)$}\label{sec:strat}

Let $X$ be a closed hyperbolic cone-3-manifold with all cone angles smaller than $2\pi$; we denote as usual its regular part by $M$ and its meridian set by $\vec{\mu}$.
Its deformation space $C_{-1}(X)$ is defined as the space of all hyperbolic cone-manifold structures on the same underlying topological space (but possibly with different singular loci) and with the same regular part $M$.
Let $C_{-1}(X)_{comp}$ be the union $\bigcup C_{-1}(X,\vec{\nu})$ over all the families $\vec{\nu}$ such that $\vec{\mu}\cup\vec{\nu}$ is an admissible pair-of-pants decomposition
(note that  $C_{-1}(X,\vec{\nu}) =  C_{-1}(X,\vec{\nu}')$ if each splittable curve of $\vec{\nu}$ is equivalent to a splittable curve of $\vec{\nu}'$ and reciprocally).
This space of compatible deformations contains a large piece of the neighbourhood of $X$ in $C_{-1}(X)$.
Theorems \ref{thm:construction} and \ref{thm:main} give a description of the shape of $C_{-1}(X)_{comp}$, and its expression as $\bigcup C_{-1}(X,\vec{\nu})$ shows that it is a stratified space near $X$: 

\begin{itemize}
\item The (closed) top-dimensional strata have dimension $3g-3$ and correspond to the $C_{-1}(X,\vec{\nu})$-spaces for which all the curves in $\vec{\nu}$ are splittable and disjoint. In these strata, we have the local parametrization $(\alpha,\ell): C_{-1}(X, \vec{\nu}) \rightarrow (0,2\pi)^N \times \R_{\geq 0}^{3g-3-N}$ given by Theorem \ref{thm:main}. 

\item We can associate a lower-dimensional stratum to each family of curves $\vec{\nu}' \in \coprod N_j$ such that $\vec{\mu}\cup\vec{\nu}'$ is a subset of an admissible pair-of-pants decomposition of $\partial \bar M_\varepsilon$, and such that the curves in $\vec{\nu}'$ are splittable and disjoint.
This (closed) stratum then consists of the deformations that are only split along the curves of $\vec{\nu}'$, and it also admits a parametrization by the cone angles of the original edges and the lengths of the new ones.
An alternate description is as the intersection $\bigcap C_{-1}(X,\vec{\nu})$ over all the family of curves $\vec{\nu}$ containing $\vec{\nu}'$ and such that $\vec{\mu}\cup\vec{\nu}$ is admissible;
note however that not all of the $C_{-1}(X,\vec{\nu})$-spaces in this intersection correspond to top-dimensional strata.
Indeed, it is even possible that $\vec{\nu}'$ cannot be completed by any disjoint, non-homotopic splittable curve, in which case the corresponding stratum does not lie in any higher-dimensional one.

\item At the intersection of all these $C_{-1}(X,\vec{\nu})$-spaces we find of course the bottom-dimensional stratum $C_{-1}(X,\Sigma)$, consisting of the deformations that do not split any vertex. 
\end{itemize}

Locally, $C_{-1}(X)$ can contain other types of deformation, besides the obvious extension of the above description to include generalized splittable curves. First of all, as we have seen in section \ref{sec:split}, it may be possible to split $X$ along curves with trivial holonomy. However, the parametrization results are no longer valid in that case. Secondly, a curve $\gamma$ on a link's smooth part $N_j$ may be ``weakly  splittable'', if $N_j$ lies on the boundary of the open set of spherical structures for which $\gamma$ is splittable (cf.~Lemma \ref{lem:open}). In that case, even though $X$ cannot be directly split along $\gamma$, it can happen that arbitrarily small deformations of it admit such splittings. The resulting elements of $C_{-1}(X)$ then do not form a well-defined stratum. 
The fifth picture on Fig.~\ref{fig:Lissajous} provides an example of a weakly splittable curve: under small perturbations of the spherical structure, it can become splittable; the resulting cone angle is either close to $3\pi$ or to $5\pi$, and these two possibilities correspond to different signs for the direction of the Hamiltonian flow, cf.~Lemma \ref{lem:sign}.

\subsubsection*{A complex of curves realization}

The local stratified structure on $C_{-1}(X)_{comp}$ is in fact closely related to the geometry of the complex of curves of $\coprod N_j$. We recall briefly its definition \cite{Har,Mas-Min}: the complex of curves $K(S)$ of a surface $S$ is the simplicial complex whose $n$-dimensional simplices correspond to sets of $n+1$ homotopy classes of non-peripheral simple closed curves on $S$, realizable as disjoint curves. In particular, its top-dimensional simplices (or facets) are in bijection with the pair-of-pants decompositions of $S$. Note that the inclusion of $\coprod N_j$ in $\partial \bar M_\varepsilon$ induces an embedding of $K(\coprod N_j)$ into $K(\partial \bar M_\varepsilon)$ as the \emph{link} of the simplex $K_{\vec{\mu}}$ corresponding to the family of meridians of $X$. 

Every stratum of $C_{-1}(X)_{comp}$ is determined by a family of curves $\vec{\nu}'$, and thus corresponds to a simplex $K_{\vec{\nu}'}$ in $K(\coprod N_j)$. But this correspondence is not one-to-one: we have seen that there can exist on $\coprod N_j$ two homotopic, non-equivalent splittable curves; splitting deformations along one curve and along the other belong to two disjoint strata of $C_{-1}(X)_{comp}$, that both correspond to the same simplex of the curve complex. Besides, most families of disjoint homotopy classes on $\coprod N_j$ cannot be realized by splittable curves.

However, we can still use the results of this article to construct a geometric realization of a modified curve complex.
Let $K'$ be the subcomplex of $K(\coprod N_j)$ obtained by removing all faces whose family of homotopy classes cuts $\coprod N_j$ into subsurfaces at least one of which has reducible holonomy representation, and also all faces whose family of homotopy classes contains one with trivial holonomy (since it obviously cannot  be completed into an admissible pair-of-pants decomposition). Otherwise,  Proposition \ref{prop:admdecomp} asserts that any remaining simplex is included in a facet, i.e.~a $(3g-4-N)$-dimensional one. Note however that generically $K'=K(\coprod N_j)$ and this restriction is unnecessary. We consider the cone $C(K')$ over $K'$; it is not a simplicial complex, but its truncated version carries this structure. Now let us choose an arbitrary function $s$ which associates to each homotopy class $[\gamma]$ on $\coprod N_j$ a sign $s([\gamma]) \in \{-1,1\}$. We will see that any choice of such a function yields a continuous map $\Xi_s$ from a neighbourhood  $\mathcal{V}$ of $(\alpha,0)$ in $(0,2\pi)^N\times C(K')$ to $\Defo(M) \stackrel{loc}{\simeq} X(\pi_1M,\SL_2(\C))$, where $\alpha \in (0,2\pi)^N$ is the vector of the cone angles of $X$.

Every top-dimensional simplex of $K'$ determines (up to homotopy) a family of curves $\vec{\nu}$ such that $\mathcal{C}=\vec{\mu}\cup \vec{\nu}$ is an admissible pair-of-pants decomposition of $\partial \bar M_\varepsilon$. This simplex $K'_{\vec{\nu}}$ can be parametrized as $\{(x_1,\ldots,x_{3g-3-N}) \in \R_{\geq 0}^{3g-3-N}\ :\ \sum_j x_j=1\}$, where the subset $\{x_i=1\}$ is the vertex associated to the single curve $\nu_i$. By Theorem \ref{thm:param}, $\mathcal{C}$ gives rise to a local coordinate chart $\Phi_\mathcal{C}$ on $X(\pi_1M,\SL_2(\C))$, and as in Proposition \ref{prop:param3}, we can consider the local diffeomorphism $\Psi_{\vec{\nu},s} : \Phi_{\mathcal{C}}^{-1}((-2,2)^N\times \R_s^{3g-3-N}\times\{0\}^{3g-3-N}) \to (0,2\pi)^N \times \R_{\geq 0}^{3g-3-N}$,
\begin{multline*}
\Psi_{\vec{\nu},s}(\chi) = \Big(f_1^{-1}(\tr_{\mu_1}(\chi)),\ldots,f_N^{-1}(\tr_{\mu_N}(\chi)),
 s([\nu_1]) \Im \tau_{\nu_1}(\chi) \sqrt{4-\tr_{\nu_1}(\chi)^2},\\ \ldots,s([\nu_{3g-3-N}]) \Im \tau_{\nu_{3g-3-N}}(\chi) \sqrt{4-\tr_{\nu_{3g-3-N}}(\chi)^2}\Big).
\end{multline*}
Let $\mathcal{V}_{\vec \nu} = \mathcal{V} \cap \left((0,2\pi)^N \times C(K'_{\vec \nu})\right)$. We can now define the continuous map $\Xi_s$ 
by requiring that its restriction to $\mathcal{V}_{\vec \nu}$ is given by
\begin{eqnarray*}
\Xi_{\vec \nu,s} : \mathcal{V}_{\vec \nu} \subset (0,2\pi)^N \times K'_{\vec \nu} \times (0,\infty) &\to &X(\pi_1M,\SL_2(\C)) \\
(\alpha'_1,\ldots,\alpha'_N,x_1,\ldots,x_{3g-3-N},t) &\mapsto & \Psi_{\vec{\nu},s}^{-1}(\alpha'_1,\ldots,\alpha'_N,tx_1,\ldots,tx_{3g-3-N})
\end{eqnarray*}
It is then easy to check that these restrictions coincide on $\mathcal{V}_{\vec \nu}\cap \mathcal{V}_{\vec{ \nu'}}$ as soon as it is non-empty. We remark that this map cannot be proper, because $K'$ is locally infinite whereas $\Defo(M)$ is locally compact.

There is an obvious relation between this map $\Xi_s$ and the stratified structure on $C_{-1}(X)_{comp}$. Indeed, the holonomy characters of the elements of $C_{-1}(X,\vec{\nu})$ correspond to the image by $\Xi_s$ of $(0,2\pi)^N \times C(K'_{\vec \nu})$, as soon as the function $s$ maps the curves in $\vec{\nu}$ to their signs as defined in Lemma \ref{lem:sign}. The non-properness of $\Xi_s$ then corresponds to the existence of converging sequences of representations, as discussed below.

\subsubsection*{Limit points}

It is remarkable that the stratified structure of $C_{-1}(X)_{comp}$ is \emph{not locally finite}. Indeed, we have seen in section \ref{sec:examples} an example with infinitely many non-homotopic splittable curves. It implies that we can have infinitely many different strata intersecting at a lower-dimensional one. However, $C_{-1}(X)$ is immersed in $\Defo(M)$, which is locally compact; this means that in this situation, there must exist a sequence of holonomy characters converging to a limit. 

\begin{figure}
\centering 
\includegraphics{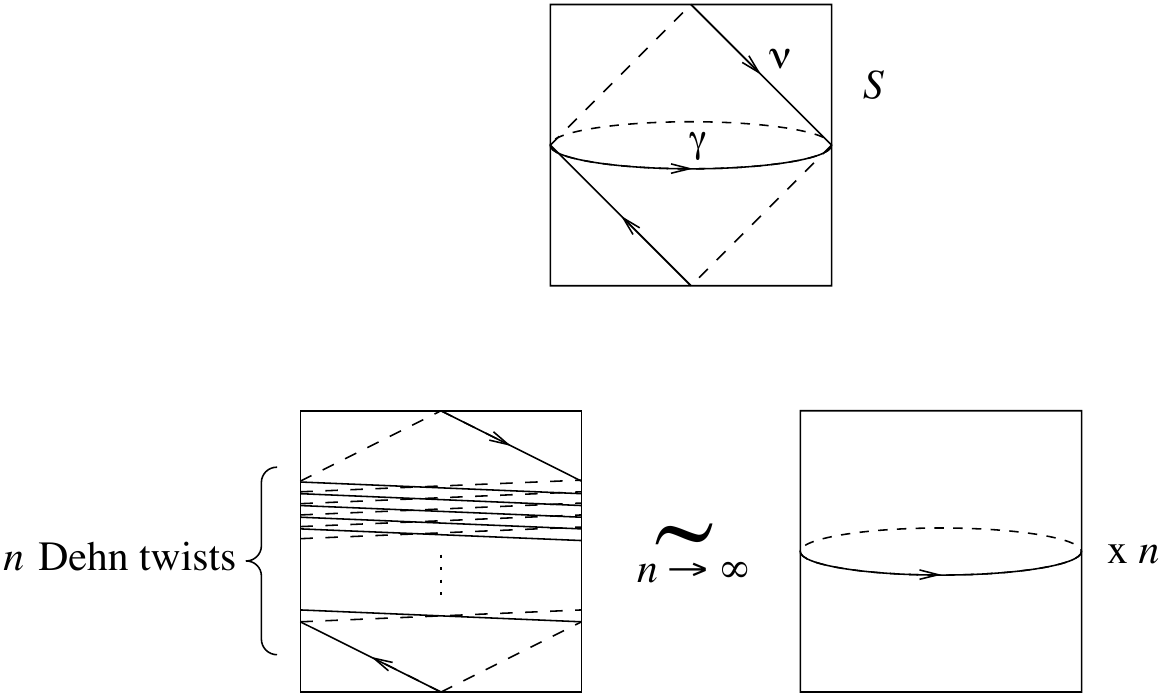}
\caption{Converging deformations}\label{fig:limit}
\end{figure}

An exact computation can be done in the example of section \ref{sec:examples}. We recall that $S$ is the double of a spherical square with angle $2\pi/3$ (see Fig.~\ref{fig:limit}). We have seen that the image $\nu_n$ of the curve $\nu$ after $n$ (full) Dehn twists along the curve $\gamma$ satisfies the splitting condition and has cone angle $2n\pi+4\pi/3$ if $n$ is even, and $2n\pi+2\pi/3$ if $n$ is odd. Let $\rho_n$ be the representation obtained by splitting $S$ along $\nu_n$ by a length of $l/n$ for some positive real number $l$. Then one can compute that this sequence of representations converges to a limit representation $\rho_\infty$, which is exactly the representation obtained by splitting $S$ along $\gamma$ by a length of $\frac{\sqrt{3}}{2}l$. This result is more striking expressed in terms of the complex Hamiltonian flow: the sequence of characters obtained from $[\rho]$ by following the flow associated to $\tr_{\nu_n}$ for a complex time $z/n$, converges to the character obtained by following the flow associated to $\tr_{\gamma}$ for a complex time $z$. This means that in some sense, after a large number of Dehn twists along $\gamma$, the curve $\nu_n$ becomes equivalent to $n$ copies of $\gamma$ itself.

\subsection{The polyhedral case}

If $P$ is a convex hyperbolic polyhedron, then its double $D(P)$ has a natural cone-manifold structure with cone angles smaller than $2\pi$, so the results of this article can be applied to $D(P)$. However for most splittings of the vertices, the resulting cone-manifold will no longer be the double of a polyhedron; this contrasts strongly with the non-splitting case, see \cite{Montcouq2}. 

For simplicity we will work with {\em face-marked} polyhedra, that is, polyhedra equipped with a bijective map from the set of faces $F$ to $\{1,2,\ldots,|F|\}$. Let $\overline{Pol}(n)$ be the set of strictly convex, face-marked, hyperbolic polyhedra with $n$ faces.  Each element of this set can be defined as a (non-redundant) intersection of $n$ half-spaces. Using the correspondence between oriented half-spaces of $\H^3$ and points of the de Sitter space $dS^3$, we see that $\overline{Pol}(n)$ can be identified with an open subset of $(dS^3)^n$; in particular, it is a smooth manifold of dimension $3n$. The isometry group $\Isom(\H^3)=  \Isom(dS^3) = \SO_0(1,3)$ acts freely and discontinuously on this set, so that the quotient $Pol(n)$, the space of congruence classes of (face-marked) strictly convex hyperbolic polyhedra with $n$ faces, is a smooth manifold of dimension $3n-6$. Using the Euler formula and the equality $2N =  \sum_{j=1}^k m_j$ (recall that $N$ is the number of edges, $k$ the number of vertices, and $m_j$ the valence of the $j$-th vertex), we obtain that the dimension of $Pol(n)$ is also equal to $N+\sum_{j=1}^k (m_j-3) = 3g-3$. This is exactly the dimension of $C_{-1}(D(P),\vec{\nu})$, for any family of curves $\vec{\nu}$ satisfying the assumptions of Theorem \ref{thm:main}.

In a neighbourhood of $D(P)$, the double construction embeds naturally $Pol(n)$ into $C_{-1}(D(P))$, but as a union of several $C_{-1}(D(P),\vec{\nu})$-spaces. Actually, it is rather simple to determine which families of curves $\vec{\nu}$ yield deformations that are double of polyhedra. The link of each vertex of $D(P)$ is the double of a convex spherical polygon, and the splitting curves must preserve this double structure. Thus $C_{-1}(D(P),\vec{\nu})$ corresponds to polyhedral deformations if and only if each curve $\nu_i$ in $\vec{\nu}$ is the double of an arc joining two non-adjacent edges of one of those spherical polygons. It is easy to see that such curves $\nu_i$ satisfy the splitting condition and give admissible pair-of-pants decompositions. Theorem \ref{thm:main} then yields, for each such $\vec{\nu}$, a parametrization by the dihedral angles and new edges' lengths of the corresponding subset of $Pol(n)$ in a neighbourhood of $D(P)$.  

\begin{figure}
\centering
\includegraphics[scale=1]{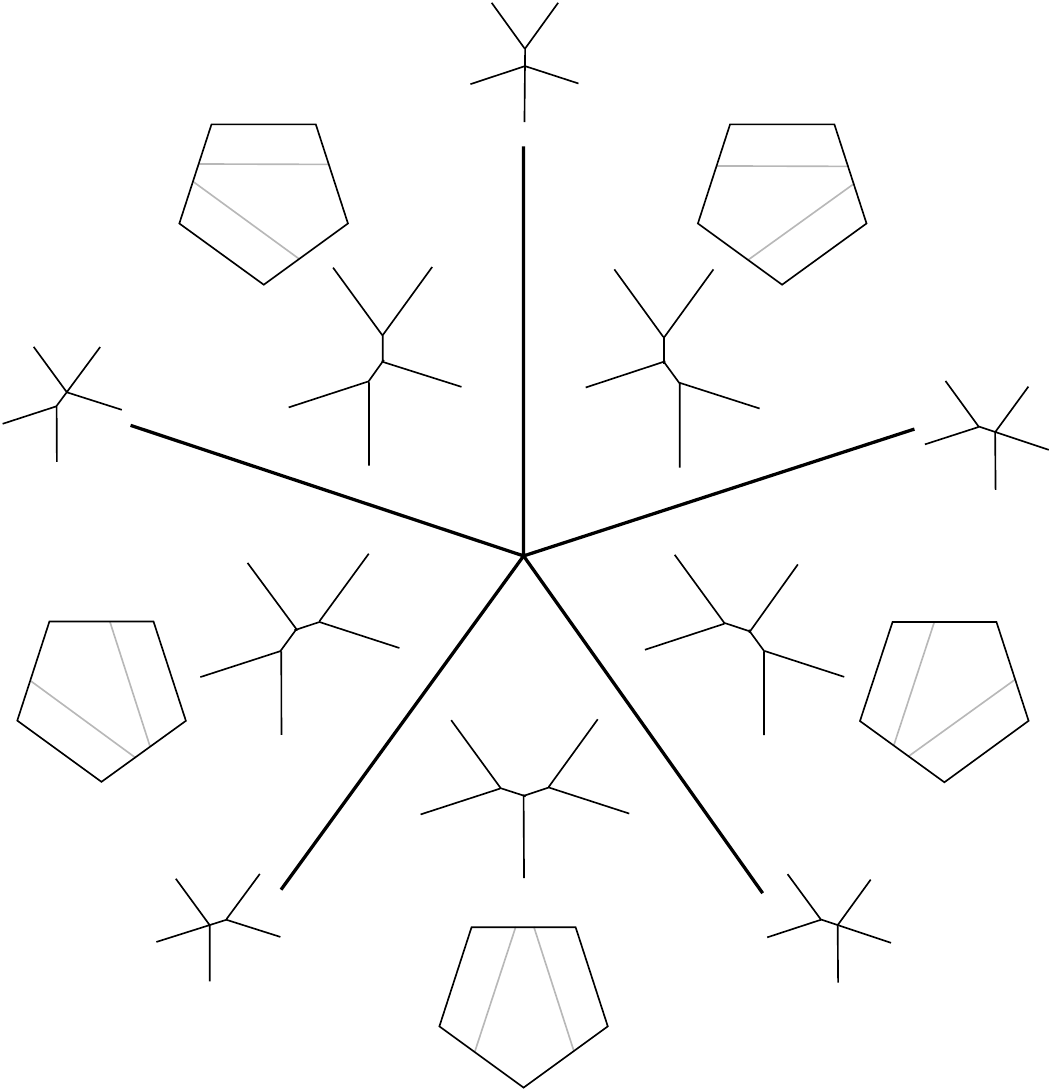}
\caption{Partition of $Pol(n)$ corresponding to all the possible splittings of a valence 5 vertex.}
\label{fig:split5}
\end{figure}

We give in Fig.~\ref{fig:split5} a schematic picture of $Pol(n)$ near $D(P)$ for a polyhedron $P$  having one vertex of valence 5 and all others of valence 3 (e.g.~a pyramid with a pentagonal base). The deformations modifying the dihedral angles are not depicted, but should be thought as an $N$-dimensional space perpendicular to the plane of the figure. Locally, $Pol(n)$ is the smooth union of five $C_{-1}(D(P),\vec{\nu})$-spaces, corresponding to the five different ways of splitting the spherical pentagon whose double is the spherical link of the valence 5 vertex. Any other pair-of-pants decomposition of this spherical link (satisfying the assumptions of Theorem \ref{thm:main}) will yield deformations of $D(P)$ that are no longer double of polyhedra.

\end{document}